\documentclass[11pt]{amsart}

\usepackage{ulem}
\usepackage{graphicx}
\usepackage{subfigure}
\usepackage{amssymb}
\usepackage[top=25mm, bottom=25mm, left=25mm, right=25mm]{geometry}
\usepackage{color}
\usepackage{bm}
\usepackage{bbm}
\usepackage{epstopdf}

\newcommand{\vct}[1]{\bm{\mathsf{#1}}}
\newcommand{\pvct}[1]{\bm{#1}}
\newcommand{\mtx}[1]{\bm{\mathsf{#1}}}

\newcommand{\pxx}{\pvct{x}}
\newcommand{\pyy}{\pvct{y}}

\newcommand{\uu}{\vct{u}}
\newcommand{\vv}{\vct{v}}

\newcommand{\TT}{\vct{T}}

\newcommand{\scalefactor}{t}
\newcommand{\ngauss}{q}

\theoremstyle{definition}
\newtheorem{remark}{Remark}

\numberwithin{definition}{section}

\newcommand{\lsp}{\vspace{3mm}}

\begin{document}

\begin{center}
\textbf{An accelerated Poisson solver based on multidomain spectral discretization}

\lsp

\textit{\small T.~Babb, A.~Gillman, S.~Hao, P.G.~Martinsson}\\
\textit{\small Department of Applied Mathematics, University of Colorado at Boulder}

\lsp

\lsp

\begin{minipage}{143mm}
\textbf{Abstract:}
This paper presents a numerical method for variable coefficient elliptic
PDEs with mostly smooth solutions on two dimensional domains.
The PDE is discretized via a multi-domain spectral collocation method
of high local order (order 30 and higher have been tested and work well).
Local mesh refinement results in highly
accurate solutions even in the presence of local irregular behavior due
to corner singularities, localized loads, etc.
The system of linear equations
attained upon discretization is solved using a direct (as opposed to iterative)
solver with $O(N^{1.5})$ complexity for the factorization stage and $O(N \log N)$
complexity for the solve.
The scheme is ideally suited for executing the elliptic solve required when
parabolic problems are discretized via time-implicit techniques.
In situations where the geometry remains unchanged between time-steps,
very fast execution speeds are obtained since the solution operator
for each implicit solve can be pre-computed.
\end{minipage}

\end{center}

\section{Introduction}
This manuscript describes a direct solver for elliptic PDEs with
variable coefficients, such as, e.g.,
\begin{equation}
\label{eq:basic}
\left\{\begin{aligned}
\mbox{}[Au](\pxx) =&\ g(\pxx),\qquad &\pxx \in \Omega,\\
   u(\pxx) =&\ f(\pxx),\qquad &\pxx \in \Gamma,
\end{aligned}\right.
\end{equation}
where $A$ is a variable coefficient elliptic differential operator
\begin{multline}
\label{eq:defA}
[Au](\pxx) = -c_{11}(\pxx)[\partial_{1}^{2}u](\pxx)
-2c_{12}(\pxx)[\partial_{1}\partial_{2}u](\pxx)
-c_{22}(\pxx)[\partial_{2}^{2}u](\pxx)\\
+c_{1}(\pxx)[\partial_{1}u](\pxx)
+c_{2}(\pxx)[\partial_{2}u](\pxx)
+c(\pxx)\,u(\pxx),
\end{multline}
where $\Omega$ is a rectangular domain in $\mathbb{R}^{2}$ with boundary $\Gamma = \partial \Omega$,
where all coefficient functions ($c$, $c_{i}$, $c_{ij}$) are smooth,
and where $f$ and $g$ are given functions. The generalization to domains that are either unions
of rectangles, or can via local parameterizations be mapped to a union of rectangles is
relatively straight-forward \cite[Sec.~6.4]{2012_spectralcomposite}.
The technique is specifically developed to accelerate implicit time stepping techniques for
parabolic PDEs such as, e.g., the heat equation
\begin{equation}
\label{eq:heat}
\left\{\begin{aligned}
\Delta u(\pxx,t) =&\ \frac{\partial u}{\partial t}(\pxx,t),\qquad &&\pxx \in \Omega,\ t > 0,\\
   u(\pxx,t) =&\ f(\pxx,t),\qquad &&\pxx \in \Gamma,\ t > 0,\\
   u(\pxx,0) =&\ g(\pxx),\qquad &&\pxx \in \Omega.
\end{aligned}\right.
\end{equation}
When (\ref{eq:heat}) is discretized using an implicit time-stepping scheme (e.g.~backwards Euler or Crank-Nicolson),
one is required to solve for each time-step an equation of the form (\ref{eq:basic}), see Section \ref{sec:timestep}.
With the ability to combine very high order discretizations with a highly efficient means of time-stepping parabolic
equations, we believe that the proposed method will be particularly well suited for numerically solving the Navier-Stokes
equation at low Reynolds numbers.

The proposed solver is direct and builds an approximation to the solution operator of
(\ref{eq:basic}) via a hierarchical divide-and-conquer approach. It is conceptually
related to classical nested dissection and multifrontal methods \cite{2006_davis_directsolverbook,1989_directbook_duff,george_1973},
but provides tight integration between the direct solver and the discretization procedure.
Locally, the scheme relies on high order spectral discretizations,
and collocation of the differential operator. We observe that while classical nested dissection
and multifrontal solvers slow down dramatically as the discretization order is increased
\cite[Table 3]{2013_martinsson_DtN_linearcomplexity}, the
proposed method retains high efficiency regardless of the discretization order.
The method is an evolution of the scheme described in
\cite{2012_spectralcomposite,2012_martinsson_composite_orig}, and
later refined in \cite{2013_martinsson_ItI,2013_martinsson_DtN_linearcomplexity,2016_martinsson_HPS_3D}.
One novelty of the present
work is that it describes how problems with body loads can be handled efficiently
(the previous papers \cite{2012_spectralcomposite,2013_martinsson_ItI,2013_martinsson_DtN_linearcomplexity}
consider the case where $g=0$ in (\ref{eq:basic})).
A second novelty is that local
mesh refinement is introduced to enable the method to accurately solve problems
involving concentrated loads, singularities at re-entrant corners, and other
phenomena that lead to localized loss of regularity in the solution.
(In contrast, the previous papers
\cite{2012_spectralcomposite,2013_martinsson_ItI,2013_martinsson_DtN_linearcomplexity,2016_martinsson_HPS_3D}
restrict attention to uniform grids.)

The principal advantage of the proposed solver, compared to commonly used solvers for (\ref{eq:basic}),
is that it is \textit{direct} (as opposed to \textit{iterative}), which makes it particularly well
suited for problems for
which efficient pre-conditioners are difficult to find, such as, e.g., problems with oscillatory solutions.
The cost to build the solution operator is in the most basic version of the scheme $O(N^{3/2})$, where
$N$ is the number of discretization points.
However, the practical efficiency of the solver is very high and the superlinear scaling is hardly visible
until $N > 10^7$. When the number of discretization points is higher than $10^7$, the scheme can be modified
to attain linear complexity by implementing techniques analogous to those described in \cite{2013_martinsson_ItI}.
Once the solution operator has been built, the time required to apply it to execute a solve given a boundary
condition and a body load is either $O(N\log N)$ for the basic scheme,
or $O(N)$ for the accelerated scheme, with a small scaling constant in either case. In
Section \ref{sec:num}, we demonstrate
that even when $N=10^6$, the time for solving (\ref{eq:basic}) with a precomputed
solution operator is approximately one second on a standard office laptop.

The discretization scheme we use is related to earlier work on spectral collocation methods
on composite (``multi-domain'') grids, such as, e.g., \cite{1998_kopriva_multidomain,2000_hesthaven_pseudospectral},
and in particular Pfeiffer \textit{et al} \cite{2003_pfeiffer_spectralmultidomain}.
The differences and similarities between the various techniques is discussed in
detail in \cite{2012_spectralcomposite}. Our procedure is also conceptually related to
so-called ``reduction to the interface'' methods, see \cite{2004_khoromskij_reduction_to_interface}
and the references therein.  Such ``interface'' methods also use local solution operators defined
on boundaries but typically rely on variational formulations of the PDE, rather than the
collocation techniques that we employ.

The manuscript is organized as follows:
Section \ref{sec:outline} provides a high level description of the proposed method.
Sections \ref{sec:leaf} and \ref{sec:merge} describe the local discretization scheme.
Section \ref{Algorithm} describes the nested dissection type solver used to solve
the system of linear equations resulting from the discretization.
Section \ref{sec:refinement} describes how local mesh refinement can be introduced to the scheme.
Section \ref{sec:num} provides results from numerical experiments that establish the
efficiency of the proposed method.

\section{Overview of algorithm}
\label{sec:outline}

\begin{figure}[b]
\textit{Level 0} \hspace{20mm}
\textit{Level 1} \hspace{20mm}
\textit{Level 2} \hspace{20mm}
\textit{Level 3} \hspace{20mm}
\textit{Level 4}

\includegraphics[width=\textwidth]{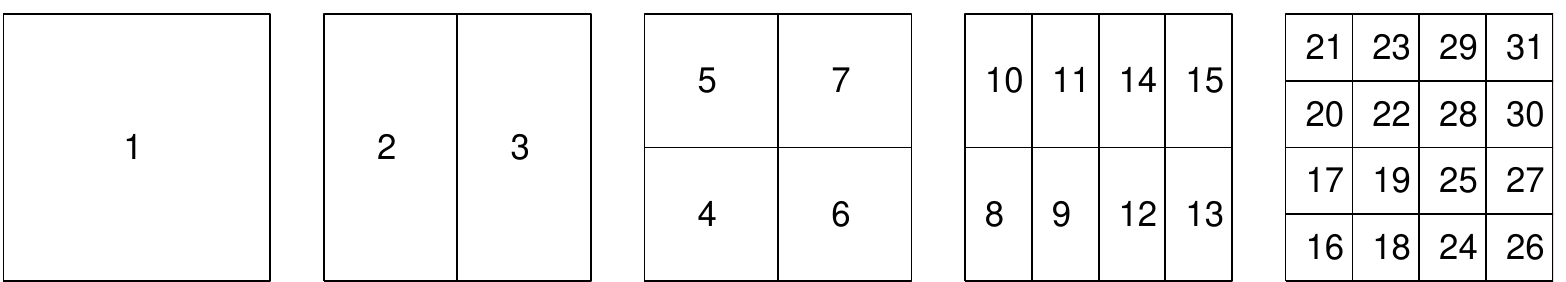}

\caption{The square domain $\Omega$ is split into $4 \times 4$ leaf boxes.
These are then gathered into a binary tree of successively larger boxes
as described in Section \ref{sec:outline}. One possible enumeration
of the boxes in the tree is shown, but note that the only restriction is
that if box $\tau$ is the parent of box $\sigma$, then $\tau < \sigma$.}
\label{fig:tree_numbering}
\end{figure}

The proposed method is based on a hierarchical subdivision of the computational domain,
as illustrated in Figure \ref{fig:tree_numbering} for the case of $\Omega = [0,1]^2$.
In the uniform mesh version of the solver, the tree of boxes is built by recursively splitting the original box in halves.
The splitting continues until each box is small enough that the solution, and its
first and second derivatives, can accurately be resolved on a local tensor product grid
of $p\times p$ Chebyshev nodes  (where, say, $p=10$ or $p=20$).

Once the tree of boxes has been constructed, the actual solver consists of two stages.
The first, or ``build'',  stage consists of a single upwards pass through the tree of
boxes, starting with the leaves and going up to larger boxes. On each leaf, we place
a local $p\times p$ tensor product grid of Chebyshev nodes, and then discretize the
restriction of (\ref{eq:basic}) via a basic collocation scheme, as in \cite{2000_trefethen_spectral_matlab}.
By performing dense linear algebraic operations on matrices of size at most $p^{2} \times p^{2}$,
we form for each leaf a local solution operator and an approximation to the local
Dirichlet-to-Neumann (DtN) operator, as described in Section \ref{sec:leaf}.
The build stage then continues with an upwards pass through the tree
(going from smaller boxes to larger) where for each parent box, we construct approximations
to its local solution operator and its local DtN operator by ``merging'' the corresponding
operators for its children, cf.~Section \ref{sec:merge}. The end result of the ``build stage''
is a hierarchical representation
of the overall solution operator for (\ref{eq:basic}). Once this solution operator is available,
the ``solve stage'' takes as input a given boundary data $f$ and a body load $g$, and constructs
an approximation to the solution $u$ valid throughout the domain via two passes through the
tree: first an upwards pass (going from smaller boxes to larger) where ``particular solutions''
that satisfy the inhomogeneous equation are built, and then a downwards pass where the boundary
conditions are corrected.

The global grid of collocation points used in the upwards and downwards passes is obtained
by placing on the edge of each leaf a set of $q$ Gaussian interpolation nodes (a.k.a.~Legendre
nodes). Observe that this parameter $q$ is in principle distinct from the local parameter $p$
which specifies the order of the local Chebyshev grids used to construct the solution operators
on the leaves. However, we typically choose $p=q+1$ or $p=q+2$.


\section{Leaf computation}
\label{sec:leaf}

In this section, we describe how to numerically build the various
linear operators (represented as dense matrices) needed for a given leaf $\Omega_{\tau}$ in the
hierarchical tree. To be precise, let $u$ be the solution to the local equation
\begin{equation}
\label{eq:local}
\left\{\begin{aligned}
\mbox{}[Au](\pxx) =&\ g(\pxx),\qquad &\pxx \in \Omega_{\tau},\\
         u(\pxx)  =&\ d(\pxx),\qquad &\pxx \in \Gamma_{\tau},
\end{aligned}\right.
\end{equation}
for some given (local) Dirichlet data $d$. We then build approximations
to two linear operators that both take $d$ and $g$ as their inputs. The
first operator outputs the local solution $u$ on $\Omega_{\tau}$ and the
second outputs the boundary fluxes of $u$ on $\Gamma_{\tau}$.

\subsection{Notation}
We work with
two sets of interpolation nodes on the domain $\Omega_\tau$. First, let
$\{\pyy_{j}\}_{j=1}^{4q}$ denote the nodes obtained by placing
$q$ Gaussian nodes on each of the four sides of $\Omega_{\tau}$.
Next, let $\{\pxx_{i}\}_{i=1}^{p^2}$ denote the nodes in a
$p\times p$ Chebyshev grid on $\Omega_{\tau}$.
We partition the index vector for the nodes in the Chebyshev grid as
$$
\{1,2,\dots,p^{2}\} = I_{\rm ce} \cup I_{\rm ci}
$$
so that $I_{\rm ce}$ holds the (\underline{C}hebyshev)
\underline{e}xterior nodes and $I_{\rm ci}$ holds the
(\underline{C}hebyshev) \underline{i}nterior nodes.
Let $\uu_{\rm c}$, $\uu_{\rm ci}$, $\uu_{\rm ce}$, and $\uu_{\rm ge}$
denote vectors holding approximations to the values of the solution
$u$ at the interpolation nodes:
$$
\uu_{\rm c} \approx \{u(\pxx_{i})\}_{i=1}^{p^{2}},\qquad
\uu_{\rm ci} \approx \{u(\pxx_{i})\}_{i\in I_{\rm ci}},\qquad
\uu_{\rm ce} \approx \{u(\pxx_{i})\}_{i\in I_{\rm ce}},\qquad
\uu_{\rm ge} \approx \{u(\pyy_{j})\}_{j=1}^{4q}.
$$
Let $\vv_{\rm ge} \in \mathbb{R}^{4q}$ denote a vector holding
boundary fluxes of $u$ on the Gaussian grid, so that
\begin{align*}
\vv_{\rm ge}(j) &\approx [\partial_{1}u](\pyy_{hj})\qquad\mbox{when}\ \pyy_{j}\ \mbox{lies on a vertical boundary},\\
\vv_{\rm ge}(j) &\approx [\partial_{2}u](\pyy_{hj})\qquad\mbox{when}\ \pyy_{j}\ \mbox{lies on a horizontal boundary}.
\end{align*}
Observe that our sign convention for boundary fluxes means that a positive flux sometimes represents
flow into the box and sometimes out of the box.
Finally, let $\vct{d}_{\rm ge}$ and $\vct{g}_{\rm ci}$ denote tabulations
of the boundary data and the body load,
$$
\vct{d}_{\rm ge} = \{d(\pyy_{j})\}_{j}^{4q},\qquad
\vct{g}_{\rm ci} = \{g(\pxx_{i})\}_{i\in I_{\rm ci}}.
$$
Our objective is now to construct the matrices that map $\{\vct{d}_{\rm ge},\vct{g}_{\rm ci}\} $ to
$\vv_{\rm ge}$ and $\uu_{\rm c}$.

\subsection{Discretization on the Cheyshev grid}
In order to execute the local solve on $\Omega_{\tau}$ of (\ref{eq:local}),
we use a classical spectral collocation technique, as described, e.g., in
\cite{2000_trefethen_spectral_matlab}. To this end, let $\mtx{D}^{(1)}$
and $\mtx{D}^{(2)}$ denote the $p^2 \times p^2$ spectral differentiation
matrices on the $p\times p$ Chebyshev grid. (In other words, for any
function $u$ that is a tensor product of polynomials of degree at most
$p-1$, the differentiation matrix \textit{exactly} maps a vector of
collocated function values to the vector of collocated values of its
derivative.) Further, let $\mtx{A}$ denote the matrix
$$
\mtx{A} =
- \mtx{C}_{11}(\mtx{D}^{(1)})^2 - 2\mtx{C}_{12}\mtx{D}^{(1)}\mtx{D}^{(2)} - \mtx{C}_{22}(\mtx{D}^{(2)})^2
+ \mtx{C}_{1}\mtx{D}^{(1)} + \mtx{C}_{2}\mtx{D}^{(2)} + \mtx{C},
$$
where $\mtx{C}_{ij}$ are diagonal matrices with entries $\{c_{ij}(\pxx_k)\}_{k=1}^{p^{2}}$,
and $\mtx{C}_{i}$ and $\mtx{C}$ are defined analogously.
Next, partition the matrix $\mtx{A}$ to separate interior and exterior nodes via
$$
\mtx{A}_{\rm ci,ci} = \mtx{A}(I_{\rm ci},I_{\rm ci}),	
\qquad\mbox{and}\qquad
\mtx{A}_{\rm ci,ce} = \mtx{A}(I_{\rm ci},I_{\rm ce}).
$$
Collocating (\ref{eq:local}) at the interior nodes then results in the discretized equation
\begin{equation}
\mtx{A}_{\rm ci,ci}\,\vct{u}_{\rm ci} +
\mtx{A}_{\rm ci,ce}\,\vct{d}_{\rm ce} =
\vct{g}_{\rm ci},
\label{eq:local_sol}
\end{equation}
where $\vct{d}_{\rm ce} = \{d(\pvct{x}_{i})\}_{i \in I_{\rm ce}}$ encodes the local
Dirichlet data $d$.

\subsection{Solving on the Chebyshev grid}
\label{sec:localsolve}
While solving (\ref{eq:local_sol}) gives the solution at the interior Chebychev nodes, it does not
give a map to the boundary fluxes $\vct{v}_{\rm ge}$ that we seek.  These are found by following
the classic approach of writing the solution as the superposition of the homogeneous and
particular solutions.  Specifically, the solution to (\ref{eq:local}) is split as
$$
u = w + \phi
$$
where $w$ is a \textit{particular solution}
\begin{equation}
\label{eq:part}
\left\{\begin{aligned}
Aw(\pxx) =&\ g(\pxx),\qquad &\pxx \in \Omega_{\tau},\\
 w(\pxx) =&\ 0,\qquad &\pxx \in \Gamma_{\tau},
\end{aligned}\right.
\end{equation}
and where $\phi$ is a \textit{homogeneous solution}
\begin{equation}
\label{eq:hom}
\left\{\begin{aligned}
A\phi(\pxx) =&\ 0,\qquad &\pxx \in \Omega_{\tau},\\
 \phi(\pxx) =&\ d(\pxx),\qquad &\pxx \in \Gamma_{\tau}.
\end{aligned}\right.
\end{equation}
Discretizing (\ref{eq:part}) on the Chebyshev grid, and collocating at the internal nodes, we get the equation
$$
\mtx{A}_{\rm ci,ce}\vct{w}_{\rm ce} + \mtx{A}_{\rm ci,ci}\vct{w}_{\rm ci} = \vct{g}_{\rm ci}.
$$
Observing that $\vct{w}_{\rm ce} = 0$, the particular solution is given by
\begin{equation}
\label{eq:alpine1}
\vct{w}_{\rm c} =
\left[\begin{array}{c} \vct{w}_{\rm ce} \\ \vct{w}_{\rm ci} \end{array}\right] =
\mtx{F}_{\rm c,ci}\vct{g}_{\rm ci},
\qquad\mbox{where}\qquad
\mtx{F}_{\rm c,ci} =
\left[\begin{array}{c} \mtx{0} \\ \mtx{A}_{\rm ci,ci}^{-1} \end{array}\right].
\end{equation}
Analogously, the discretization of (\ref{eq:hom}) on the Chebyshev grid yields
$$
\mtx{A}_{\rm ci,ce}\vct{\phi}_{\rm ce} + \mtx{A}_{\rm ci,ci}\vct{\phi}_{\rm ci} = \vct{0}.
$$
Since $\vct{\phi}_{\rm ce} = \vct{d}_{\rm ce}$, the homogeneous solution is given by
\begin{equation}
\label{eq:alpine2}
\vct{\phi}_{\rm c} =
\left[\begin{array}{c} \vct{\phi}_{\rm ce} \\ \vct{\phi}_{\rm ci} \end{array}\right] =
\left[\begin{array}{c} \mtx{I} \\ -\mtx{A}_{\rm ci,ci}^{-1}\mtx{A}_{\rm ci,ce} \end{array}\right]\vct{d}_{\rm ce}.
\end{equation}

\subsection{Interpolation and differentiation}
\label{sec:interp}

Section \ref{sec:localsolve} describes how to locally solve the BVP (\ref{eq:local}) on
the Chebyshev grid via the superposition of the homogeneous and particular solutions. Note that this
computation assumes that the local Dirichlet data $d$
is given on the \textit{Chebyshev} exterior nodes. In reality, this data will be provided
on the Gaussian nodes, and we therefore need to introduce an interpolation operator that
moves data between the different grids. To be precise, let $\mtx{L}_{\rm ce,ge}$ denote a
matrix of size $4(p-1)\times 4q$ that maps a given data vector $\vct{d}_{\rm ge}$ to
a different vector
\begin{equation}
\label{eq:alpine2p5}
\begin{array}{cccccccccc}
\vct{d}_{\rm ce} &=& \mtx{L}_{\rm ce,ge} & \vct{d}_{\rm ge}\\
4(p-1)\times 1 && 4(p-1) \times 4q & 4q \times 1
\end{array}
\end{equation}
as follows: An entry of $\vct{d}_{\rm ce}$ corresponding to an interior node is defined simply
via a standard interpolation from the Gaussian to the Chebyshev nodes on the local edge alone.
An entry of $\vct{d}_{\rm ce}$ corresponding to a corner node is defined as the average value
of the two extrapolated values from the Gaussian nodes on the two edges connecting to the corner.
(Observe that except for the four rows corresponding to the corner nodes, the matrix
$\mtx{L}_{\rm ce,ge}$ is a $4\times 4$ block diagonal matrix.)

Combining (\ref{eq:alpine1}), (\ref{eq:alpine2}), and (\ref{eq:alpine2p5}), the solution
to (\ref{eq:local}) on the Chebyshev grid is given by
\begin{equation}
\label{eq:alpine3}
\vct{u}_{\rm c} =
\vct{w}_{\rm c} + \vct{\phi}_{\rm c} =
\mtx{F}_{\rm c,ci}\,\vct{g}_{\rm ci} + \mtx{S}_{\rm c,ge}\,\vct{d}_{\rm ge},
\qquad\mbox{where}\qquad
\mtx{S}_{\rm c,ge} :=
\left[\begin{array}{c} \mtx{I}_{\rm ce,ce} \\ -\mtx{A}_{\rm ci,ci}^{-1}\mtx{A}_{\rm ci,ce} \end{array}\right]\,\mtx{L}_{\rm ce,ge}
\end{equation}
All that remains is now to determine the vector $\vct{v}_{\rm ge}$ of boundary fluxes on the
Gaussian nodes. To this end, let us define a combined interpolation and differentiation matrix
$\mtx{D}_{\rm ge,c}$ of size $4q \times p^{2}$ via
$$
\mtx{D}_{\rm ge,c} = \left[\begin{array}{c}
\mtx{L}_{\rm loc}\,\mtx{D}_{2}(I_{\rm s},:) \\
\mtx{L}_{\rm loc}\,\mtx{D}_{1}(I_{\rm e},:) \\
\mtx{L}_{\rm loc}\,\mtx{D}_{2}(I_{\rm n},:) \\
\mtx{L}_{\rm loc}\,\mtx{D}_{1}(I_{\rm w},:)
\end{array}\right],
$$
where $\mtx{L}_{\rm loc}$ is a $q\times p$ interpolation matrix from a
set of $p$ Chebyshev nodes to a set of $q$ Gaussian nodes, and where
$I_{\rm s},\,I_{\rm e},\,I_{\rm n},\,I_{\rm w}$ are four index sets, each of length $p$,
that point to the south, east, north, and west sides of the exterior nodes in the
Chebyshev grid. By differentiating the local solution on the Chebyshev grid
defined by (\ref{eq:alpine3}), the boundary fluxes $\vct{v}_{\rm ge}$ are given by
\begin{equation}
\label{eq:alpine4}
\vct{v}_{\rm ge} =
\mtx{H}_{\rm ge,ci}\,\vct{g}_{\rm ci} +
\mtx{T}_{\rm ge,ge}\,\vct{d}_{\rm ge},
\qquad\mbox{where}\qquad
\mtx{H}_{\rm ge,ci} = \mtx{D}_{\rm ge,c}\,\mtx{F}_{\rm c,ci}
\ \mbox{and}\
\mtx{T}_{\rm ge,ge} = \mtx{D}_{\rm ge,c}\,\mtx{S}_{\rm c,ge}.
\end{equation}

\section{Merging two leaves}
\label{sec:merge}

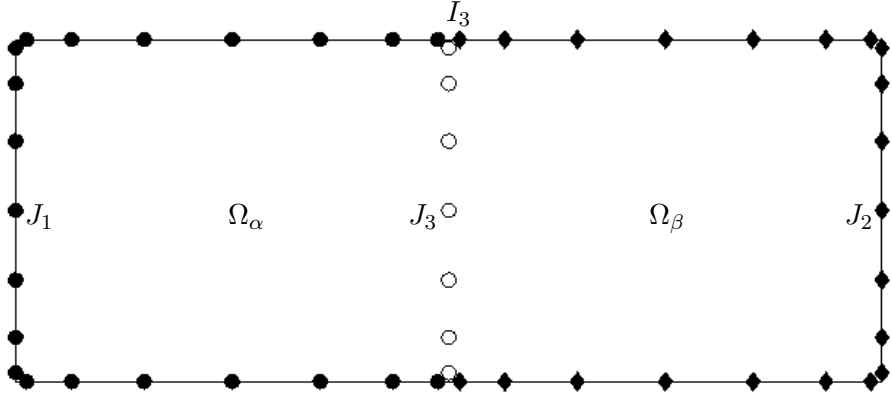
\begin{figure}
\setlength{\unitlength}{1mm}
\begin{picture}(95,55)
\put(-15,00){\includegraphics[height=55mm]{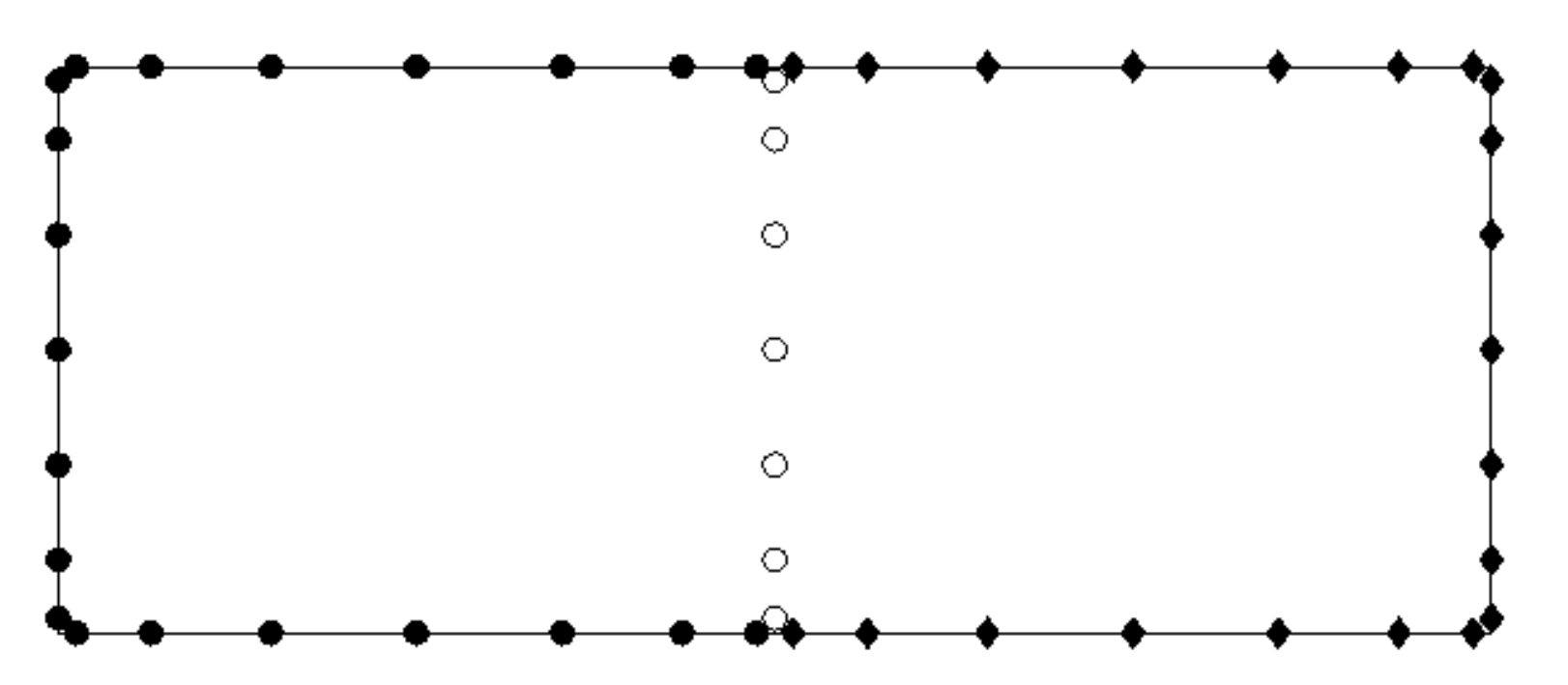}}
\put(18,25){$\Omega_{\alpha}$}
\put(74,25){$\Omega_{\beta}$}
\put(-9,25){$J_{1}$}
\put(100,25){$J_{2}$}
\put(42,25){$J_{3}$}
 \put(47,52){$I_{3}$}
\end{picture}
\caption{Notation for the merge operation described in Section \ref{sec:merge}.
Given two leaf boxes $\Omega_{\alpha}$ and $\Omega_{\beta}$, their union is denoted
$\Omega_{\tau} = \Omega_{\alpha} \cup \Omega_{\beta}$.
The sets $J_{1}$ (black circles) and $J_{2}$ (black diamonds) form the exterior nodes, while
$J_{3}$ (white circles) consists of the interior nodes.}
\label{fig:siblings_notation}
\end{figure}

Consider a rectangular box $\tau$ consisting of two leaf boxes $\alpha$ and $\beta$,
and suppose that all local operators for $\alpha$ and $\beta$ defined in Section \ref{sec:leaf}
have been computed. Our objective is now to construct the Dirichlet-to-Neumann operator
for the bigger box $\tau$ from the local operators for its children. In this operation,
only sets of Gaussian nodes on the boundaries will take part, cf.~Figure \ref{fig:siblings_notation}.
We group these nodes into three sets, indexed by vectors $J_{1}$, $J_{2}$, and $J_{3}$, defined as follows:
\begin{center}
\begin{tabular}{ll}
$J_{1}$\hspace{10mm}& Edge nodes of box $\alpha$ that are not shared with box $\beta$.\\
$J_{2}$\hspace{10mm}& Edge nodes of box $\beta$  that are not shared with box $\alpha$.\\
$J_{3}$\hspace{10mm}& Edge notes that line the interior edge shared by $\alpha$ and $\beta$.
\end{tabular}
\end{center}
We also define
$$
J_{\rm ge}^{\tau} = J_{1} \cup J_{2}
\qquad\mbox{and}\qquad
J_{\rm gi}^{\tau} = J_{3}
$$
as the exterior and interior nodes for the parent box $\tau$. Finally, we let
$\vct{h}^{\alpha}, \vct{h}^{\beta} \in \mathbb{R}^{4q}$ denote two vectors that
hold the boundary fluxes for the two local particular solutions $w^{\alpha}$ and $w^{\beta}$,
cf.~(\ref{eq:alpine4}),
\begin{equation}
\label{eq:defhleaf}
\vct{h}^{\alpha}_{\rm ge} = \mtx{H}_{\rm ge,ci}^{\alpha}\,\vct{g}_{\rm ci}^{\alpha},
\qquad\mbox{and}\qquad
\vct{h}^{\beta}_{\rm ge}  = \mtx{H}_{\rm ge,ci}^{\beta }\,\vct{g}_{\rm ci}^{\beta }.
\end{equation}
Then the equilibrium equations for each of the two leaves can be written
\begin{equation}
\label{eq:pair}
\vct{v}^{\alpha}_{\rm ge} = \mtx{T}^{\alpha}_{\rm ge,ge}\,\vct{u}^{\alpha}_{\rm ge} + \vct{h}^{\alpha}_{\rm ge},
\qquad\mbox{and}\qquad
\vct{v}^{\beta}_{\rm ge}  = \mtx{T}^{\beta }_{\rm ge,ge}\,\vct{u}^{\beta }_{\rm ge} + \vct{h}^{\beta }_{\rm ge}.
\end{equation}
Now partition the two equations in (\ref{eq:pair}) using the notation shown in
Figure \ref{fig:siblings_notation} so that
\begin{align}
\label{eq:bitter1}
\left[\begin{array}{c}
\vv_{1}\\ \vv_{3}
\end{array}\right]
=&\
\left[\begin{array}{ccc}
\mtx{T}_{1,1}^{\alpha} & \mtx{T}_{1,3}^{\alpha} \\
\mtx{T}_{3,1}^{\alpha} & \mtx{T}_{3,3}^{\alpha}
\end{array}\right]\,
\left[\begin{array}{c}
\uu_{1}\\ \uu_{3}
\end{array}\right]
+
\left[\begin{array}{c} \vct{h}_{1}^{\alpha} \\ \vct{h}_{3}^{\alpha} \end{array}\right],
\\
\label{eq:bitter2}
\left[\begin{array}{c}
\vv_{2}\\ \vv_{3}
\end{array}\right] =&\
\left[\begin{array}{ccc}
\mtx{T}_{2,2}^{\beta} & \mtx{T}_{2,3}^{\beta} \\
\mtx{T}_{3,2}^{\beta} & \mtx{T}_{3,3}^{\beta}
\end{array}\right]\,
\left[\begin{array}{c}
\uu_{2}\\ \uu_{3}
\end{array}\right]
+
\left[\begin{array}{c} \vct{h}_{2}^{\beta} \\ \vct{h}_{3}^{\beta} \end{array}\right].
\end{align}
(The subscript ``ge'' is suppressed in (\ref{eq:bitter1}) and (\ref{eq:bitter2}) since all nodes involved are Gaussian exterior nodes.)
Combine the two equations for $\vv_{3}$ in (\ref{eq:bitter1}) and (\ref{eq:bitter2}) to obtain the equation
$$
\mtx{T}_{3,1}^{\alpha}\,\uu_{1} + \mtx{T}_{3,3}^{\alpha}\,\uu_{3} + \vct{h}_{3}^{\alpha}
=
\mtx{T}_{3,2}^{\beta }\,\uu_{2} + \mtx{T}_{3,3}^{\beta }\,\uu_{3} + \vct{h}_{3}^{\beta}.
$$
This gives
\begin{equation}
\label{eq:swiss}
\uu_{3} =
\bigl(\TT^{\alpha}_{3,3} - \TT^{\beta}_{3,3}\bigr)^{-1}
\bigl( \TT^{\beta }_{3,2}\uu_{2}
      -\TT^{\alpha}_{3,1}\uu_{1}
      +\vct{h}_{3}^{\beta}
      -\vct{h}_{3}^{\alpha}
      \bigr)
\end{equation}
Using the relation (\ref{eq:swiss}) in combination with (\ref{eq:bitter1}), we find that
\begin{align*}
\left[\begin{array}{c} \vv_{1} \\ \vv_{2} \end{array}\right] =&\
\left(\left[\begin{array}{ccc}
\TT_{1,1}^{\alpha} & \mtx{0} \\
\mtx{0} & \TT_{2,2}^{\beta }
\end{array}\right] +
\left[\begin{array}{c}
\TT_{1,3}^{\alpha} \\
\TT_{2,3}^{\beta}
\end{array}\right]\,
\bigl(\TT^{\alpha}_{3,3} - \TT^{\beta}_{3,3}\bigr)^{-1}
\bigl[-\TT^{\alpha}_{3,1}\ \big|\ \TT^{\beta}_{3,2}].
\right)
\left[\begin{array}{c} \uu_{1} \\ \uu_{2} \end{array}\right] +\\
&\
\left[\begin{array}{c} \vct{h}^{\alpha}_{1} \\ \vct{h}^{\beta}_{2} \end{array}\right] +
\left[\begin{array}{c}
\TT_{1,3}^{\alpha} \\
\TT_{2,3}^{\beta}
\end{array}\right]\,
\bigl(\TT^{\alpha}_{3,3} - \TT^{\beta}_{3,3}\bigr)^{-1}
\bigl(\vct{h}_{3}^{\beta}-\vct{h}_{3}^{\alpha}\bigr).
\end{align*}
We now define the operators
\begin{align*}
\mtx{X}^{\tau}_{\rm gi,gi} =&\ \bigl(\TT^{\alpha}_{3,3} - \TT^{\beta}_{3,3}\bigr)^{-1},\\
\mtx{S}^{\tau}_{\rm gi,ge} =&\
\bigl(\TT^{\alpha}_{3,3} - \TT^{\beta}_{3,3}\bigr)^{-1}\bigl[-\TT^{\alpha}_{3,1}\ \big|\ \TT^{\beta}_{3,2}] =
\mtx{X}^{\tau}_{\rm gi,gi}\bigl[-\TT^{\alpha}_{3,1}\ \big|\ \TT^{\beta}_{3,2}],\\
\mtx{T}^{\tau}_{\rm ge,ge} =&
\left[\begin{array}{ccc}
\TT_{1,1}^{\alpha} & \mtx{0} \\
\mtx{0} & \TT_{2,2}^{\beta }
\end{array}\right] +
                    \left[\begin{array}{c}
                          \TT_{1,3}^{\alpha} \\
                          \TT_{2,3}^{\beta}
                          \end{array}\right]\,\bigl(\TT^{\alpha}_{3,3} - \TT^{\beta}_{3,3}\bigr)^{-1}\bigl[-\TT^{\alpha}_{3,1}\ \big|\ \TT^{\beta}_{3,2}] \\
 = &
\left[\begin{array}{ccc}
                          \TT_{1,1}^{\alpha} & \mtx{0}\\
                          \mtx{0} & \TT_{2,2}^{\beta }
                          \end{array}\right] +
                    \left[\begin{array}{c}
                          \TT_{1,3}^{\alpha} \\
                          \TT_{2,3}^{\beta}
                          \end{array}\right]\,\mtx{S}_{\rm gi,ge}^{\tau}.
\end{align*}


Constructing the approximate solution on the shared edge $\vct{u}_{\rm gi}^{\tau}$
can be viewed an upward pass to compute
the approximate boundary flux by

\begin{equation}
\label{eq:defhparent}
\vct{h}_{\rm ge}^{\tau} =\ \left[\begin{array}{c} \vct{h}^{\alpha}_{1} \\ \vct{h}^{\beta}_{2} \end{array}\right] +\left[\begin{array}{c}
\TT_{1,3}^{\alpha} \\
\TT_{2,3}^{\beta}
\end{array}\right]\vct{w}_{\rm gi}^{\tau},\end{equation}
where $\vct{w}_{\rm gi}^{\tau} =\ \mtx{X}_{\rm gi,gi}^{\tau}\bigl(\vct{h}_{3}^{\beta}-\vct{h}_{3}^{\alpha}\bigr)$,
followed by a downward pass
$$
\vct{u}^{\tau}_{\rm gi} = \mtx{S}_{\rm gi,ge}\vct{u}^{\tau}_{\rm ge} + \vct{w}_{\rm gi}^{\tau}.
$$

\begin{remark}[Physical interpretation of merge]
The quantities $\vct{w}^{\tau}_{\rm gi}$ and $\vct{h}_{\rm ge}^{\tau}$ have a simple
physical meaning. The vector $\vct{w}_{\rm gi}^{\tau}$ introduced above is simply a
tabulation of the particular solution $w^{\tau}$ associated with $\tau$ on the interior
boundary $\Gamma_{3}$, and $\vct{h}_{\rm ge}^{\tau}$ is the normal derivative of $w^{\tau}$.
To be precise, $w^{\tau}$ is the solution to the inhomogeneous problem, cf.~(\ref{eq:part})
\begin{equation}
\label{eq:part_tau}
\left\{\begin{aligned}
Aw^{\tau}(\pxx) =&\ g(\pxx),\qquad &\pxx \in \Omega_{\tau},\\
 w^{\tau}(\pxx) =&\ 0,\qquad &\pxx \in \Gamma_{\tau}.
\end{aligned}\right.
\end{equation}
We can re-derive the formula for $w|_{\Gamma_{3}}$ using the original mathematical operators
as follows: First observe that for $\pxx \in \Omega^{\alpha}$, we have
$A(w^{\tau} - w^{\alpha}) = g - g = 0$, so the DtN operator $T^{\alpha}$ applies to the
function $w^{\tau} - w^{\alpha}$:
$$
T^{\alpha}_{31}(w_{1}^{\tau} - w_{1}^{\alpha}) +
T^{\alpha}_{33}(w_{3}^{\tau} - w_{3}^{\alpha}) =
(\partial_{n}w^{\tau})|_{3} - (\partial_{n}w^{\alpha})|_{3}
$$
Use that $w_{1}^{\tau} = w_{1}^{\alpha} = w_{3}^{\alpha} = 0$, and that $(\partial_{n}w^{\alpha})|_{3} = h_{3}^{\alpha}$ to get
\begin{equation}
\label{eq:legomine1}
T^{\alpha}_{33}w_{3}^{\tau} = (\partial_{n}w^{\tau})|_{3} - h_{3}^{\alpha}.
\end{equation}
Analogously, we get
\begin{equation}
\label{eq:legomine2}
T^{\beta}_{33}w_{3}^{\tau} = (\partial_{n}w^{\tau})|_{3} - h_{3}^{\beta}.
\end{equation}
Combine (\ref{eq:legomine1}) and (\ref{eq:legomine2}) to eliminate $(\partial_{n}w^{\tau})|_{3}$ and obtain
$$
\bigl(T^{\alpha}_{33} - T^{\beta}_{33}\bigr) w_{3}^{\tau} =
 - h_{3}^{\alpha} + h_{3}^{\beta}.
$$
Observe that in effect, we can write the particular solution $w^{\tau}$ as
$$
w^{\tau}(\pxx) = \left\{\begin{array}{ll}
w^{\alpha}(\pxx) + \hat{w}^{\tau}(\pxx)\quad& \pxx \in \Omega^{\alpha},\\
w^{\beta }(\pxx) + \hat{w}^{\tau}(\pxx)\quad& \pxx \in \Omega^{\beta},
\end{array}\right.
$$
The function $w^{\tau}$ must of course be smooth across $\Gamma_{3}$, so the
function $\hat{w}^{\tau}$ must have a jump that exactly offsets the discrepancy
in the derivatives of $w^{\alpha}$ and $w^{\beta}$. This jump is precisely of
size $h^{\alpha} - h^{\beta}$.
\end{remark}

\section{The full solver for a uniform grid}
\label{Algorithm}

\subsection{Notation}
\label{Notation}
Suppose that we are given a rectangular domain $\Omega$, which has hierarchically
been split into a binary tree of successively smaller patches, as described
in Section \ref{sec:outline}. We then define two sets of interpolation nodes.
First, $\{\pxx_{i}\}_{i=1}^{M}$ denotes the set of nodes obtained by
placing a $p\times p$ tensor product grid of Chebyshev nodes on each leaf in
the tree. For a leaf $\tau$, let $I_{\rm c}^{\tau}$ denote an index vector
pointing to the nodes in $\{\pxx_{i}\}_{i=1}^{M}$ that lie on leaf $\tau$.
Thus the index vector for the set of nodes in $\tau$ can be partitioned into
exterior and interior nodes as follows
$$
I_{\rm c}^{\tau} = I_{\rm ce}^{\tau} \cup I_{\rm ci}^{\tau}.
$$
The second set of interpolation nodes $\{\pyy_{j}\}_{j=1}^{N}$ is obtained
by placing a set of $q$ Gaussian (``Legendre'') interpolation nodes on the
edge of each leaf. For a node $\tau$ in the tree (either a leaf or a parent),
let $I_{\rm ge}^{\tau}$ denote an index vector that marks all Gaussian
nodes that lie on the boundary of $\Omega_{\tau}$. For a parent node $\tau$,
let $I_{\rm gi}^{\tau}$ denote the Gaussian nodes that are interior to
$\tau$, but exterior to its two children (as in Section \ref{sec:merge}).

Once the collocation points have been set up, we introduce a vector
$\vct{u} \in \mathbb{R}^{M}$ holding approximations to the values of the
potential $u$ on the Gaussian collocation points,
$$
\vct{u}(j) \approx u(\pyy_{j}),\qquad j = 1,\,2,\,3,\,\dots,\,M.
$$
We refer to subsets of this vector using the short-hand
$$
\vct{u}_{\rm ge}^{\tau} = \vct{u}(I_{\rm ge}^{\tau}),
\qquad\mbox{and}\qquad
\vct{u}_{\rm gi}^{\tau} = \vct{u}(I_{\rm gi}^{\tau})
$$
for the exterior and interior nodes respectively.
At the very end of the algorithm, approximations to $u$ on the local Chebyshev
tensor product grids are constructed.
For a leaf node $\tau$, let the vectors $\vct{u}_{\rm c}^{\tau}$,  $\vct{u}_{\rm ce}^{\tau}$, and $\vct{u}_{\rm ci}^{\tau}$
denote the vectors holding approximations to the potential on sets of collocation
points in the Chebyshev grid marked by $I_{\rm c}^{\tau}$, $I_{\rm ce}^{\tau}$, and
$I_{\rm ci}^{\tau}$, respectively. Observe that these vectors are \textit{not}
subvectors of $\vct{u}$.

Before proceeding to the description of the algorithm, we introduce two sets
of auxiliary vectors. First, for any parent node $\tau$, let the vector
$\vct{w}_{\rm gi}^{\tau}$ denote the computed values of the local particular
solution $w^{\tau}$ that solves (\ref{eq:part}) on $\Omega_{\tau}$, as tabulated
on the interior line marked by $I_{\rm gi}^{\tau}$. Also, define $\vct{h}^{\tau}$
as the approximate boundary fluxes of $w^{\tau}$ as defined by (\ref{eq:defhleaf}) for a leaf
and by (\ref{eq:defhparent}) for a parent.

\subsection{The build stage}
\label{sec:build}
Once the domain is partitioned into a hierarchical tree,
we execute a ``build stage'' in which the following matrices are constructed
for each box $\tau$:

\vspace{2mm}

\begin{itemize}
\item[$\mtx{S}^{\tau}$]
For a box $\tau$, the solution operator that maps Dirichlet data $\psi$ on $\partial \Omega_{\tau}$ to values of $u$ at the interior nodes. In other
words, $\mtx{u}^{\tau}_{\rm c} = \mtx{S}^{\tau}_{\rm c,ge}\mtx{\psi}^{\tau}_{\rm ge}$ on a leaf or $\mtx{u}^{\tau}_{\rm gi} = \mtx{S}^{\tau}_{\rm gi,ge}\mtx{\psi}^{\tau}_{\rm ge}$ on a parent box.

\vspace{2mm}

\item[$\TT^{\tau}$] For a box $\tau$, the matrix that maps Dirichlet data $\psi$ on $\partial \Omega_{\tau}$ to the flux $v$ on the boundary.
In other words, $\vv^{\tau}_{\rm ge} = \TT^{\tau}_{\rm ge,ge}\mtx{\psi}_{\rm ge}$.

\vspace{2mm}

\item[$\mtx{F}^{\tau}$] For a leaf box, the matrix that maps the body load to the particular solution on the interior of the leaf assuming the Dirichlet data is zero on the boundary.  In other words $\mtx{w}^{\tau}_{\rm c} = \mtx{F}^{\tau}_{\rm c,ci}\mtx{g}_{\rm ci}$.

\vspace{2mm}

\item[$\mtx{H}^{\tau}$] For a leaf box, the matrix that maps the body load to the flux on the boundary of the leaf.  In other words $\mtx{h}^{\tau}_{\rm ge} = \mtx{H}^{\tau}_{\rm ge,ci}\mtx{g}^{\tau}_{\rm ci}$.

\vspace{2mm}

\item[$\mtx{X}^{\tau}$] For a parent box $\tau$ with children $\alpha$ and $\beta$, the matrix that maps the fluxes of the particular solution for the children on the interior of a parent to the particular solution on the interior nodes.  In other words $\mtx{w}_{\rm gi} = \mtx{X}^{\tau}_{\rm gi,gi}(\mtx{h}_{3}^{\beta}- \mtx{h}_{3}^{\alpha})$.

\end{itemize}

\vspace{2mm}

The build stage consists of a single sweep over all nodes in the tree.
Any ordering of the boxes in which a parent box is processed after its
children can be used. For each leaf box $\tau$, approximations
$\mtx{S}^{\tau}$ and $\mtx{F}^{\tau}$ to the solution operators for the
homogeneous and particular solutions are constructed.  Additionally,
approximations $\mtx{T}^{\tau}$ and $\mtx{H}^{\tau}$
to the local DtN map $\mtx{T}^{\tau}$ for the homogeneous and particular
solutions are constructed using the procedure
described in Section \ref{sec:leaf}. For a parent box $\tau$ with children
$\alpha$ and $\beta$, we construct the
solution operators $\mtx{X}_{\rm gi,gi}^{\tau}$ and
$\mtx{S}_{\rm gi,ge}^{\tau}$, and the DtN operator
$\mtx{T}_{\rm ge,ge}^{\tau}$ via the process described in Section \ref{sec:merge}.
Algorithm 1 summarizes the build stage.

\subsection{The solve stage}
\label{sec:solve}
After the ``build stage'' described in Algorithm 1 has been completed,
an approximation to the global solution operator of (\ref{eq:basic}) has been
computed, and represented through the various matrices ($\mtx{H}^{\tau}$, $\mtx{F}^{\tau}$, etc.)
described in Section \ref{sec:build}. Then given specific boundary data $f$ and
a body load $g$, the corresponding solution $u$ to (\ref{eq:basic}) can be found
through a ``solve stage'' that involves two passes through the tree, first an
upwards pass (from smaller to larger boxes), and then a downwards pass.
In the upward pass, the particular solutions and normal derivatives of the particular solution
are computed and stored in the vectors $\vct{w}$ and $\vct{h}$ respectively.  Then by
sweeping down the tree applying the solution operators $\mtx{S}$ to the Dirichlet boundary data
for each box $\tau$ and adding the particular solution, the approximate solution $\vct{u}$ is
computed. Algorithm 2 summarizes the solve stage.

We observe that the vectors $\vct{w}_{\rm gi}^{\tau}$ can all be stored on a global vector
$\vct{w} \in \mathbb{R}^{N}$. Since each boundary collocation node $\vct{y}_{j}$ belongs
to precisely one index vector $I_{\rm gi}^{\tau}$, we simply find that
$\vct{w}_{\rm gi}^{\tau} = \vct{w}(I_{\rm gi}^{\tau})$.

\begin{remark}[Efficient storage of particular solutions]
For notational simplicity, we describe Algorithm 2 (the ``solve stage'') in a way that assumes
that for each box $\tau$, we explicitly store a corresponding vector $\vct{h}_{\rm ge}^{\tau}$
that represents the boundary fluxes for the local particular solution. In practice, these vectors
can all be stored on a global vector $\vct{h} \in \mathbb{R}^{N}$, in a manner similar to how
we store $\vct{w}$. For any box $\tau$ with children $\alpha$ and $\beta$, we store on $\vct{h}$
the \textit{difference} between the boundary fluxes, so that $\vct{h}(I_{\rm gi}^{\tau}) =
-\vct{h}_{3}^{\alpha} + \vct{h}_{3}^{\beta}$. In other words, as soon as the boundary fluxes
have been computed for a box $\alpha$, we add its contributions to the vector $\vct{h}(I_{\rm ge}^{\alpha})$
with the appropriate signs and then delete it. This becomes notationally less clear, but is actually
simpler to code.
\end{remark}

%

\begin{figure}
\fbox{
\begin{minipage}{148mm}
\begin{center}
\textsc{Algorithm 1} (Build stage for problems with body load)
\end{center}

\vspace{2mm}

\begin{minipage}{146mm}
This algorithm builds all solution operators required to solve the non-homogeneous
BVP (\ref{eq:basic}).
It is assumed that if node $\tau$ is a parent of node $\sigma$, then $\tau < \sigma$.
\end{minipage}

\rule{\textwidth}{0.5pt}

\begin{tabbing}
\mbox{}\hspace{6mm} \= \mbox{}\hspace{6mm} \= \mbox{}\hspace{6mm} \= \mbox{}\hspace{6mm} \= \kill
\textbf{for} $\tau = N_{\rm boxes},\,N_{\rm boxes}-1,\,N_{\rm boxes}-2,\,\dots,\,1$\\
\> \textbf{if} ($\tau$ is a leaf)\\[1mm]
\> \> $\mtx{F}_{\rm c,ci}^{\tau} = \left[\begin{array}{c}\mtx{0}\\ \mtx{A}_{\rm ci,ci}^{-1}\end{array}\right]$
\quad\textit{[pot.] $\leftarrow$ [body load]}\\[1mm]
\> \> $\mtx{H}_{\rm ge,ci}^{\tau} = \mtx{D}_{\rm ge,c}\mtx{F}_{\rm c,ci}^{\tau}$
\quad\textit{[deriv.] $\leftarrow$ [body load]}\\[1mm]
\> \> $\mtx{S}_{\rm c,ge}^{\tau} = \left[\begin{array}{c}\mtx{I}\\ -\mtx{A}_{\rm ci,ci}^{-1}\mtx{A}_{\rm ci,ce}\end{array}\right]\,
\mtx{L}_{\rm ce,ge}$
\quad\textit{[pot.] $\leftarrow$ [pot.]}\\[1mm]
\> \> $\mtx{T}_{\rm ge,ge}^{\tau} = \mtx{D}_{\rm ge,c}\mtx{S}_{\rm c,ge}$
\quad\textit{[deriv.] $\leftarrow$ [pot.] (NfD operator)}\\[1mm]
\> \textbf{else}\\
\> \> Let $\alpha$ and $\beta$ be the children of $\tau$.\\
\> \> Partition $I_{\rm ge}^{\alpha}$ and $I_{\rm ge}^{\beta}$ into vectors $I_{1}$, $I_{2}$, and $I_{3}$ as shown in Figure \ref{fig:siblings_notation}.\\[1mm]
\> \> $\mtx{X}_{\rm gi,gi}^{\tau} = \bigl(\TT^{\alpha}_{3,3} - \TT^{\beta}_{3,3}\bigr)^{-1}$
\quad\textit{[pot.] $\leftarrow$ [deriv.]}\\[1mm]
\> \> $\mtx{S}_{\rm gi,ge}^{\tau} = \mtx{X}_{\rm gi,gi}^{\tau}
                           \bigl[-\TT^{\alpha}_{3,1}\  \big|\
                                  \TT^{\beta}_{3,2}\bigr]$
\quad\textit{[pot.] $\leftarrow$ [pot.]}\\[1mm]
\> \> $\TT_{\rm ge,ge}^{\tau} = \left[\begin{array}{ccc}
                          \TT_{1,1}^{\alpha} & \mtx{0}\\
                          \mtx{0} & \TT_{2,2}^{\beta }
                          \end{array}\right] +
                    \left[\begin{array}{c}
                          \TT_{1,3}^{\alpha} \\
                          \TT_{2,3}^{\beta}
                          \end{array}\right]\,\mtx{S}_{\rm gi,ge}^{\tau}$
\quad\textit{[deriv.] $\leftarrow$ [pot.] (NfD operator)}.\\[1mm]
\> \textbf{end if}\\
\textbf{end for}
\end{tabbing}
\end{minipage}}
\caption{Build stage.}
\label{fig:precomp_with_bodyload}
\end{figure}


\begin{figure}
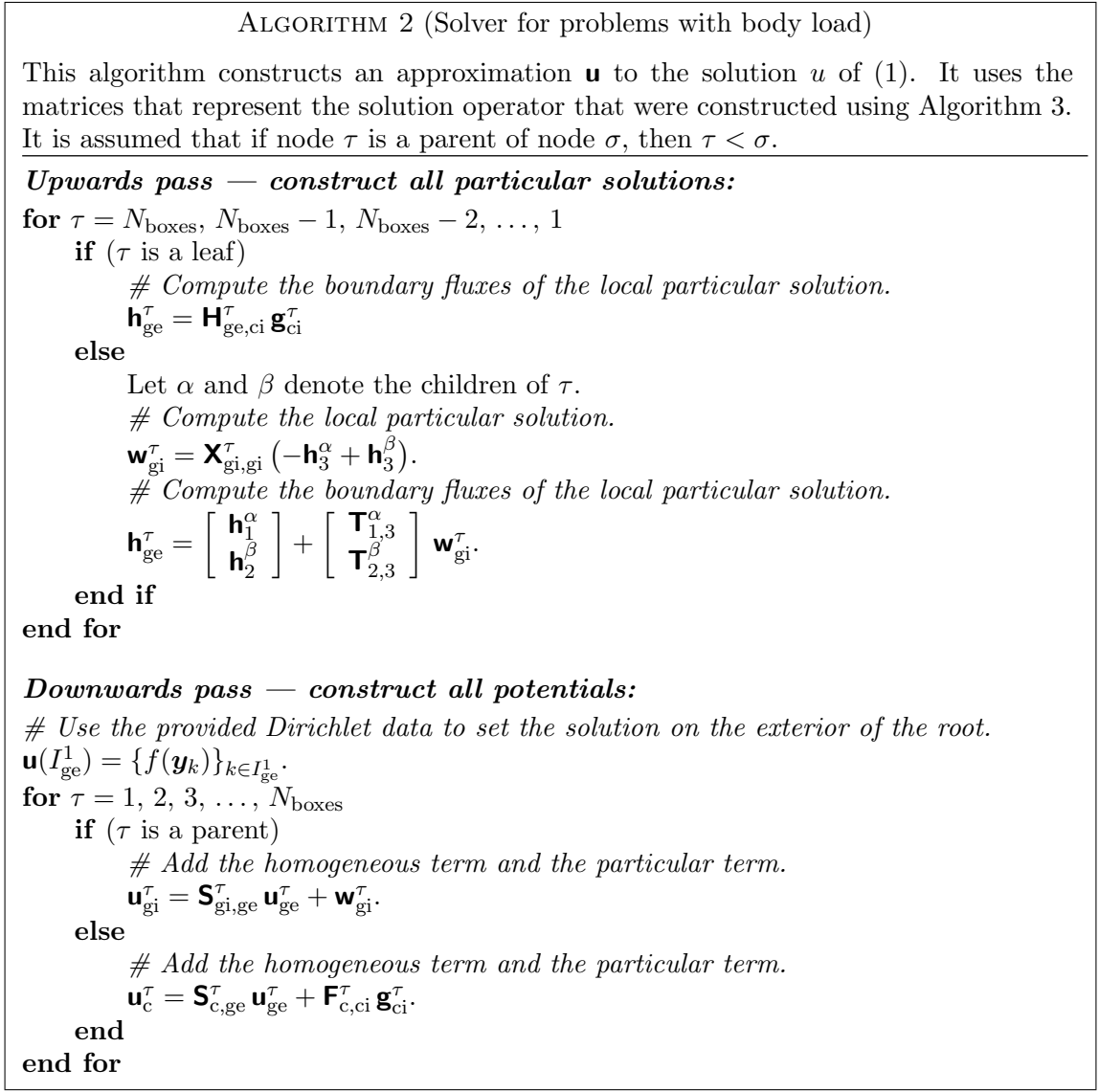

\fbox{
\begin{minipage}{148mm}
\begin{center}
\textsc{Algorithm 2} (Solver for problems with body load)
\end{center}

\vspace{2mm}

\begin{minipage}{146mm}
This algorithm constructs an approximation $\uu$ to the solution $u$ of (\ref{eq:basic}).
It uses the matrices that represent the solution operator that were
constructed using Algorithm 3.
It is assumed that if node $\tau$ is a parent of node $\sigma$, then $\tau < \sigma$.
\end{minipage}

\rule{\textwidth}{0.5pt}

\textit{\textbf{Upwards pass --- construct all particular solutions:}}
\begin{tabbing}
\mbox{}\hspace{6mm} \= \mbox{}\hspace{6mm} \= \mbox{}\hspace{6mm} \= \mbox{}\hspace{6mm} \= \kill
\textbf{for} $\tau = N_{\rm boxes},\,N_{\rm boxes}-1,\,N_{\rm boxes}-2,\,\dots,\,1$\\
\> \textbf{if} ($\tau$ is a leaf)\\
\> \> \textit{\# Compute the boundary fluxes of the local particular solution.}\\
\> \> $\vct{h}_{\rm ge}^{\tau} = \mtx{H}_{\rm ge,ci}^{\tau}\,\vct{g}_{\rm ci}^{\tau}$\\
\> \textbf{else} \\
\> \> Let $\alpha$ and $\beta$ denote the children of $\tau$.\\
\> \> \textit{\# Compute the local particular solution.}\\
\> \> $\vct{w}_{\rm gi}^{\tau} = \mtx{X}_{\rm gi,gi}^{\tau}\,\bigl(-\vct{h}_{3}^{\alpha} + \vct{h}_{3}^{\beta}\bigr).$\\
\> \> \textit{\# Compute the boundary fluxes of the local particular solution.}\\
\> \> $\vct{h}_{\rm ge}^{\tau} =
       \left[\begin{array}{c}\vct{h}_{1}^{\alpha} \\ \vct{h}_{2}^{\beta}\end{array}\right] +
       \left[\begin{array}{c}\mtx{T}_{1,3}^{\alpha} \\ \mtx{T}_{2,3}^{\beta}\end{array}\right]\,\vct{w}_{\rm gi}^{\tau}.$\\
\> \textbf{end if}\\
\textbf{end for}
\end{tabbing}

\lsp

\textit{\textbf{Downwards pass --- construct all potentials:}}
\begin{tabbing}
\mbox{}\hspace{6mm} \= \mbox{}\hspace{6mm} \= \mbox{}\hspace{6mm} \= \mbox{}\hspace{6mm} \= \kill
\textit{\# Use the provided Dirichlet data to set the solution on the exterior of the root.}\\
$\uu(I_{\rm ge}^{1}) = \{f(\pyy_{k})\}_{k \in I_{\rm ge}^{1}}$.\\
\textbf{for} $\tau = 1,\,2,\,3,\,\dots,\,N_{\rm boxes}$\\
\> \textbf{if} ($\tau$ is a parent)\\
\> \> \textit{\# Add the homogeneous term and the particular term.}\\
\> \> $\uu_{\rm gi}^{\tau} = \mtx{S}_{\rm gi,ge}^{\tau}\,\uu_{\rm ge}^{\tau} + \vct{w}_{\rm gi}^{\tau}$.\\
\> \textbf{else}\\
\> \> \textit{\# Add the homogeneous term and the particular term.}\\
\> \> $\uu_{\rm c}^{\tau} = \mtx{S}_{\rm c,ge}^{\tau}\,\uu_{\rm ge}^{\tau} +
                            \mtx{F}_{\rm c,ci}^{\tau}\,\vct{g}_{\rm ci}^{\tau}$.\\
\> \textbf{end}\\
\textbf{end for}
\end{tabbing}
\end{minipage}}
\caption{Solve stage.}
\label{fig:solver_with_bodyload}
\end{figure}

\subsection{Algorithmic complexity}
In this section, we determine the asymptotic complexity of the direct solver.  The analysis is very similar to the analysis seen in \cite{2013_martinsson_DtN_linearcomplexity} for no body load.
Let $N_{\rm leaf} = 4\ngauss$ denote the number of Gaussian nodes on the boundary
of a leaf box, and let $p^2$ denote the number of Chebychev nodes used in the leaf computation.
In the asymptotic analysis, we set $p=q+2$, so that $p \sim q$.
Let $L$ denote the number of levels in the binary tree.  This means there are
$2^L$ boxes.  Thus the total number of discretization nodes $N$ is approximately $2^L\ngauss^{2}$.

In processing a leaf, the dominant cost involves matrix inversion (or factorization followed by triangular
solve) and matrix-matrix multiplications. The largest matrices encountered are of size $O(q^2) \times O(q^2)$,
making the cost to process one leaf $O(q^6)$.  Since there are $N/q^2$ leaf boxes, the total
cost of pre-computing approximate DtN operators for all the bottom level is
$\sim (N/q^2)\times q^6\sim N\,q^4$.

Next, consider the process of merging two boxes, as described in Section \ref{sec:merge}.
On level $\ell$, there are $2^{\ell}$ boxes, that each have $O(2^{-\ell/2}N^{0.5}))$
nodes along their boundaries. (On level $\ell=2$, there are
$4$ boxes that each have side length one half of the original side length; on level
$\ell=4$, there are $16$ boxes that have side length one quarter of the original side
length; etc.) The cost of executing a merge is dominated by the cost to perform matrix
algebra (inversion, multiplication, etc) of dense matrices of size $2^{-\ell/2}N^{0.5}\times  2^{-\ell/2}N^{0.5}$.
This makes the total cost for the merges in the upwards pass
$$
\sum_{\ell=1}^{L} 2^\ell\times \left(2^{-\ell/2}N^{0.5}\right)^3 \sim
\sum_{\ell=1}^{L} 2^\ell\times 2^{-3\ell/2}N^{1.5} \sim
N^{1.5}\sum_{\ell=1}^{L} 2^{-\ell/2} \sim N^{1.5}.
$$

Finally, consider the cost of the solve stage (Algorithm 2).
We first apply at each of the $2^{L}$ leaves the operators $\mtx{H}_{\rm ge,ci}^{\tau}$,
which are all of size $4q \times (p-2)^{2}$, making the overall cost
$\sim 2^{L}q^{3} = N\,q$ since $p\sim q$ and $N\sim 2^{L}q^{2}$.
In the upwards sweep, we apply at level $\ell$ matrices of size $O(2^{-\ell/2}N^{0.5}) \times O(2^{-\ell/2}N^{0.5})$
on $2^{\ell}$ boxes, adding up to an overall cost of
$$
\sum_{\ell=1}^{L}2^{\ell}\times \left(2^{-\ell/2}\,N^{0.5}\right)^{2} \sim
\sum_{\ell=1}^{L}2^{\ell}\times 2^{-\ell}\,N \sim
\sum_{\ell=1}^{L}N \sim NL\sim N\log N.
$$
The cost of the downwards sweep is the same. However, the application of
the matrices $\mtx{F}_{\rm c,ci}^{\tau}$ at the leaves is more expensive
since these are of size $O(q^{2}) \times O(q^{2})$, which adds up to an
overall cost of $2^{L}\,q^{4} = N\,q^{2}$.



The analysis of the asymptotic storage requirements perfectly mirrors the analysis
of the flop count for the solve stage, since each matrix that is stored is used
precisely once in a matrix-vector multiplication. In consequence, the amount of
storage required is
\begin{equation}
\label{eq:mem}
R \sim N\,q^{2} + N\,\log N.
\end{equation}

\begin{remark}[A storage efficient version]
\label{remark:econ}
The storage required for all solution operators can become prohibitive when the local order
$q$ is high, due to the term $N\,q^{2}$ in (\ref{eq:mem}).
One way to handle this problem is to not store the local solution
operators for a leaf, but instead simply perform small dense solves each time the ``solve
stage'' is executed. This makes the solve stage slower, obviously, but has the benefit of
completely eliminating the $Nq^{2}$ term in (\ref{eq:mem}). In fact, in this modified version,
the overall storage required is $\sim NL \approx N\,\log_{2}(N/q^{2})$, so we see that the storage costs
decrease as $q$ increases (as should be expected since we do all leaf computations from scratch
in this case). Figure \ref{Memory} provides numerical results illustrating the memory requirements
of the various approaches.
\end{remark}

\section{Local refinement}
\label{sec:refinement}

When solving a boundary value problem like (\ref{eq:basic}) it is common to have a localized
loss of regularity due to, e.g., corners on the boundary, a locally non-smooth body load or
boundary condition, or a localized loss of regularity in the coefficient functions in the
differential operator. A common approach to efficiently regain high accuracy without excessively
increasing the number of degrees of freedom used, is to locally refine the mesh near the troublesome
location. In this manuscript, we assume the location is known and given, and that we manually
specify the degree of local refinement. The difficulty that arises is that upon refinement,
the collocation nodes on neighboring patches do not necessarily match up.  To remedy this,
interpolation operators are introduced to transfer information between patches.
(The more difficult problem of determining how to automatically detect regions
that require mesh refinement is a topic of current research.)


\subsection{Refinement criterion} 
\begin{figure}
\setlength{\unitlength}{1mm}
\begin{picture}(95,80)
\put(-30,00){\includegraphics[height=75mm]{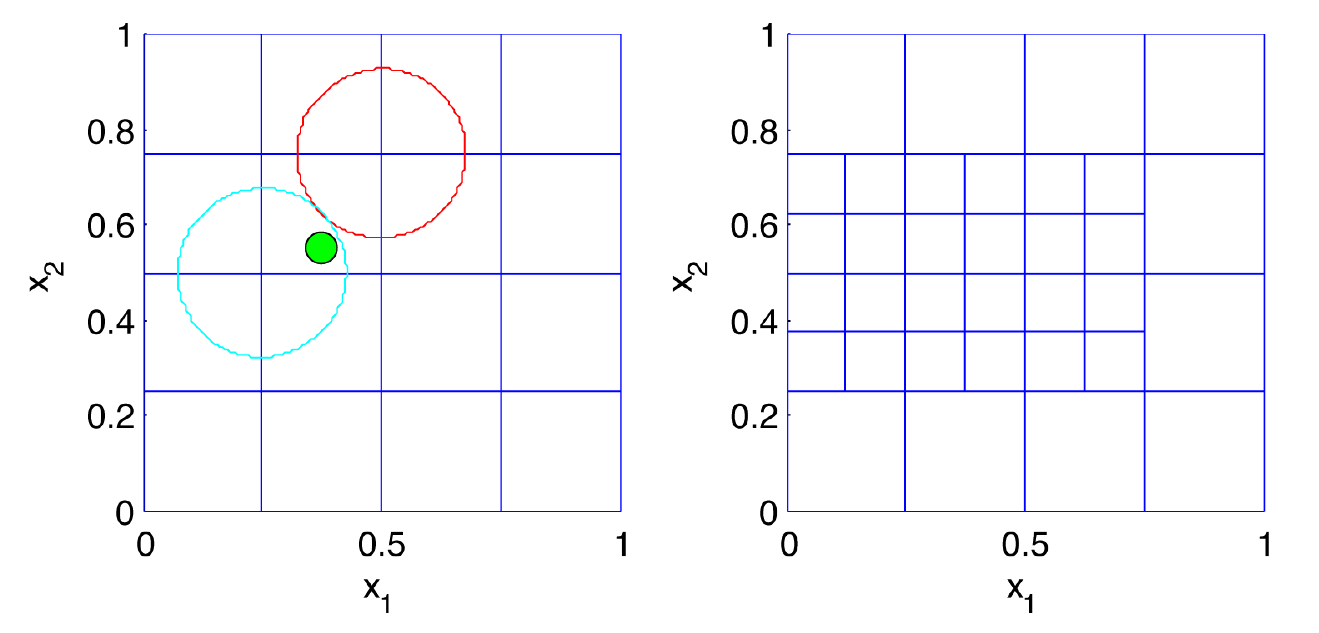}}
\put(-8,35){$\Omega_{\gamma}$}
\put(20,60){$\Omega_{\beta}$}
\end{picture}
\caption{A sample domain where we desire to refine the grid at $\hat{\pxx}$, shown by the green circle.  For leaf box $\Omega_{\gamma}$ the shortest distance to $\hat{\pxx}$ satisfies $d_{\gamma} < \scalefactor l_{\gamma}$, so $\Omega_{\gamma}$ is refined.  The maximum distance $\scalefactor l_{\gamma}$ is shown by the blue circle, which is centered at the closest point from $\Omega_{\gamma}$ to $\hat{\pxx}$.  For leaf box $\Omega_{\beta}$ the shortest distance to $\hat{\pxx}$ does not satisfy $d_{\beta} < \scalefactor l_{\beta}$, so $\Omega_{\beta}$ is not refined.  The maximum distance $\scalefactor l_{\beta}$ is shown by the red circle, which is centered at the closest point from $\Omega_{\beta}$ to $\hat{\pxx}$.}
\label{Close}
\end{figure}

Suppose we desire to refine our discretization at some point $\hat{\pvct{x}}$ in the computational domain
(the point $\hat{\pvct{x}}$ can be either in the interior or on the boundary).
Consider as an example the situation depicted in Figure \ref{Close}.  For each level of refinement, we split any leaf box that
contains $\hat{\pvct{x}}$ and any ``close'' leaf boxes into a $2 \times 2$ grid of equal-sized leaf boxes.  In Figure \ref{Close} we perform one level of refinement and find there are 6 leaf boxes ``close'' to $\hat{\pvct{x}}$, which is represented by the green dot.  These 6 boxes are refined into smaller leaf boxes.

A leaf box $\Omega_{\tau}$ is close to $\hat{\pvct{x}}$ if the distance $d_{\tau}$ from $\hat{\pvct{x}}$ to the box $\Omega_{\tau}$ satisfies
$d_{\tau} \leq \scalefactor l_{\tau}$, where $\scalefactor = \sqrt{2}$ and $2l_{\tau}$ is the length of one side of the leaf box $\Omega_{\tau}$.  In Figure \ref{Close} we show circles of size $\scalefactor l_{\gamma}$ and $\scalefactor l_{\beta}$ at the points in $\Omega_{\gamma}$ and $\Omega_{\beta}$ closest to $\hat{\pvct{x}}$ (in this case the boxes are all the same size so $l_{\gamma} = l_{\beta}$).  We see $\hat{\pvct{x}}$ is ``close'' to $\Omega_{\gamma}$, but not ``close'' to $\Omega_{\beta}$.
Just as in section \ref{Notation}, we place a $p\times p$ tensor product grid of Chebyshev nodes on each
new leaf and a set of $q$ Gaussian (``Legendre'') interpolation nodes on the
edge of each leaf. The vector $\{\pyy_{j}\}_{j=1}^{N}$ holds the locations of all Gaussian nodes across all
leaves in the domain.



\subsection{Refined mesh}
\label{RefCriteria}

\begingroup
\begin{figure}[htbp]
\subfigure[]{\includegraphics[width=0.4\textwidth]{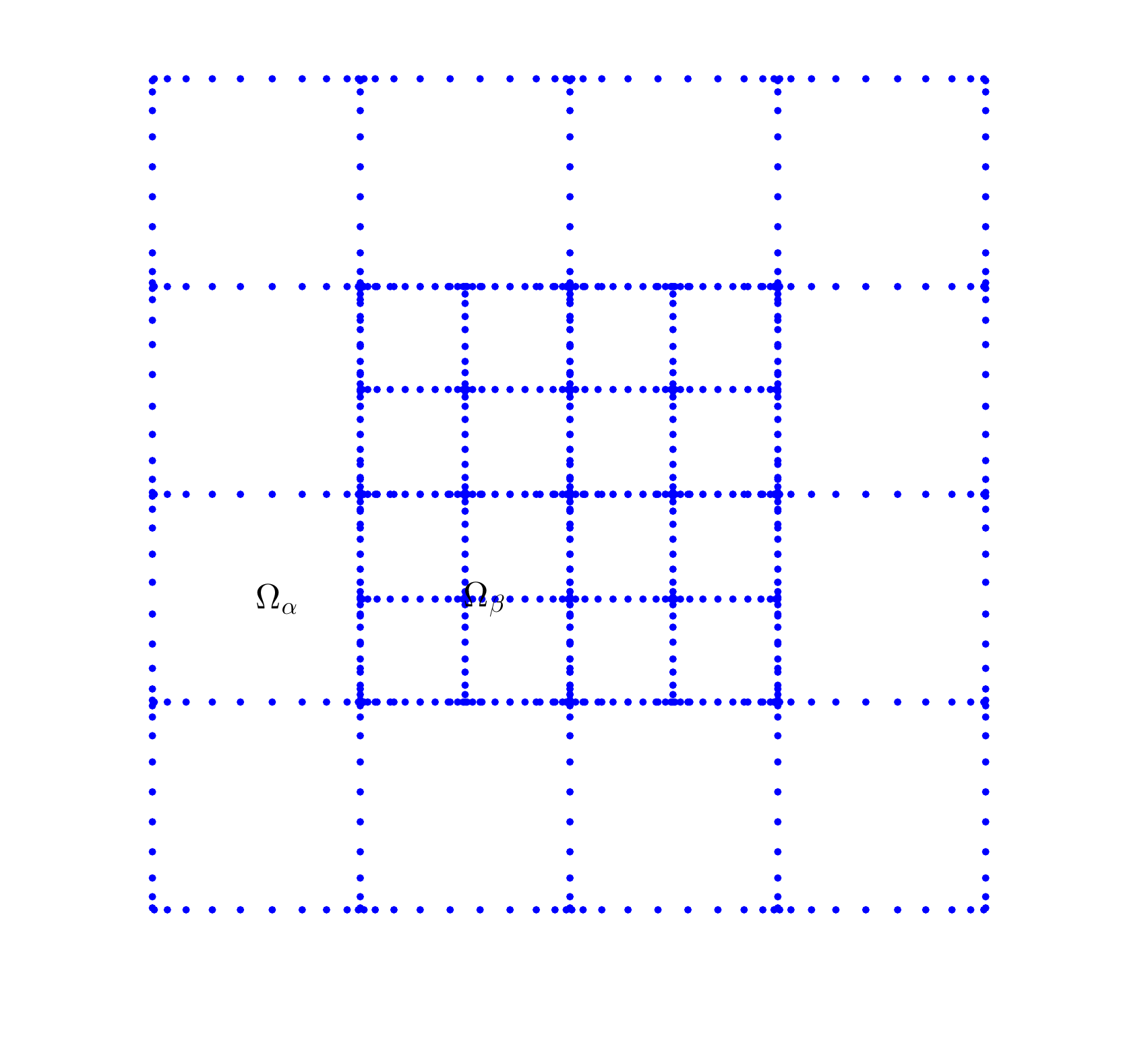}}
\subfigure[] {\includegraphics[width=0.5\textwidth]{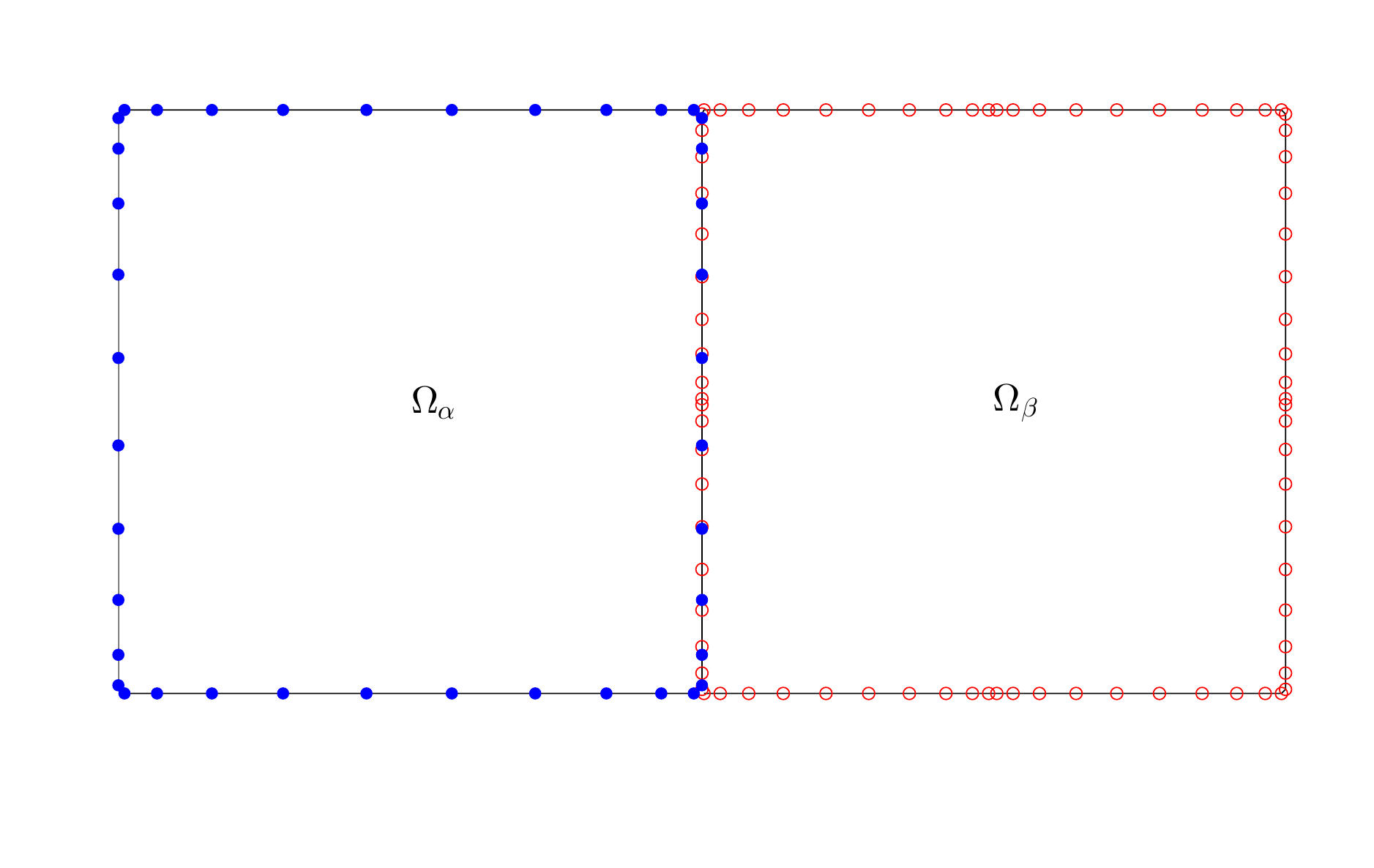}}
\label{fig:1b}
\caption{(a) Grid with refinement at the center. (b) A close up of neighbors $\Omega_{\alpha}$ and $\Omega_{\beta}$.  Since only one of the boxes is refined the exterior Gaussian nodes on the common boundary are not aligned.} \label{fig:interp_fig}
\end{figure}
\endgroup

Notice that with the refined grid the nodes along common boundaries are no longer aligned.
Figure \ref{fig:interp_fig} is an example of such a grid.   This is a
problem during the build stage of the method since the merge operation is performed by equating
the Neumann data on the common boundary.  We begin the discussion on how to address this problem
by establishing some notation.


Define two boxes as \textit{neighbors} if they are on the same level of the tree and they are adjacent.
In the case that only one of two neighbors has been refined, such as $\Omega_\alpha$ and $\Omega_\beta$ in Figure \ref{fig:interp_fig},
special attention needs to be paid to the nodes on the common boundary.  In order to merge boxes
with different number of Gaussian nodes on the common edge, interpolation operators will be
required.  The next section describes this process in detail.

Consider the nodes on the common boundary between the two leaf boxes $\Omega_\alpha$ and $\Omega_\beta$.  Let $q$ denote the number of
Gaussian nodes on one side of each leaf.  Let $\{\mtx{J}_{\alpha,i}\}_{i=1}^{q}$ denote the
index vector for the common boundary nodes from $\Omega_{\alpha}$ and $\{\mtx{J}_{\beta,i}\}_{i=1}^{2q}$
denote the index vector for the common boundary nodes from $\Omega_{\beta}$.  That is, recalling that $\mtx{y}$ holds the locations of all Gaussian nodes in the domain, $\mtx{y}(\mtx{J}_{\alpha,i})$ contains the $q$ Gaussian nodes on the Eastern side of box $\Omega_{\alpha}$ and $\mtx{y}(\mtx{J}_{\beta,i})$ contains the $2q$ nodes on the Western side of box $\Omega_{\beta}$.


\subsection{Modifications to build stage}

Once the grid with the Gaussian and Chebyshev nodes is constructed, as described in section \ref{RefCriteria},
the build stage starts with the construction of all leaf operators as described
in Section \ref{sec:leaf}.  Then, for simplicity of presentation, boxes are merged from
the lowest level moving up the tree.
After merging the children of a refined parent, such as $\Omega_{\beta}$ in Figure \ref{fig:interp_fig} it is seen that the parent's exterior nodes do not align with the exterior nodes of any neighbor which has not been refined.

Recalling the index notation used in section \ref{RefCriteria}, we form the interpolation matrix $\mtx{P}_{\rm up,W}$ mapping data on $\mtx{y}(\mtx{J}_{\alpha,i})$ to data on $\mtx{y}(\mtx{J}_{\beta,i})$ and the interpolation matrix $\mtx{P}_{\rm down,W}$ mapping data on $\mtx{y}(\mtx{J}_{\beta,i})$ to data on $\mtx{y}(\mtx{J}_{\alpha,i})$. Observe that when interpolating from two sets of $q$ Gaussian nodes to a set of $q$ Gaussian nodes, the interpolation must be done as two separate interpolations from $q$ to $q/2$ nodes.  The matrix $\mtx{P}_{\rm down,W}$ is a block diagonal matrix consisting of two $q/2 \times q$ matrices (assuming $q$ is divisible by 2).

For a refined parent, such as $\Omega_{\beta}$ in Figure \ref{fig:interp_fig}, we form the operators $\mtx{T}^{\beta}$, $\mtx{S}^{\beta}$, and $\mtx{X}^{\beta}$ and form the interpolation operators for every side of the parent, regardless of whether the exterior nodes align with the neighbor's exterior nodes.  Observe that in the case of $\Omega_{\beta}$ in Figure \ref{fig:interp_fig} the Eastern and Northern sides of $\Omega_{\beta}$ will have $\mtx{P}_{\rm down} = \mtx{P}_{\rm up}= \mtx{I}$, where $\mtx{I}$ is the identity matrix.

Then interpolation operators mapping the entire boundary data between fine and coarse grids are
given by block diagonal matrices $\mtx{P}_{\rm{up}}$ and $\mtx{P}_{\rm{down}}$ whose diagonal
blocks are the interpolation operators for each edge.  The interpolation operators for $\Omega_\beta$ are

$$\mtx{P}^{\beta}_{\rm{up}}=\texttt{blkdiag}(\mtx{P}^{\beta}_{\rm{up},S},\mtx{P}^{\beta}_{\rm{up},E},\mtx{P}^{\beta}_{\rm{up},N},\mtx{P}^{\beta}_{\rm{up},E})$$
and
$$\mtx{P}^{\beta}_{\rm{down}}=\texttt{blkdiag}(\mtx{P}^{\beta}_{\rm{down},S},\mtx{P}^{\beta}_{\rm{down},E},\mtx{P}^{\beta}_{\rm{down},N},\mtx{P}^{\beta}_{\rm{down},E}).$$
(The text  \texttt{blkdiag} denotes the function that forms a block diagonal matrix from its arguments.)
Then we form the new operators $\mtx{T}^{\beta}_{\rm new}$ and $\mtx{S}^{\beta}_{\rm new}$ for the parent box $\Omega_{\beta}$ as follows
$$\mtx{T}^{\beta}_{\rm new} = \mtx{P}^{\beta}_{\rm down} \mtx{T}^{\beta} \mtx{P}^{\beta}_{\rm up}$$
and
$$\mtx{S}^{\beta}_{\rm new} = \mtx{S}^{\beta} \mtx{P}^{\beta}_{\rm up}.$$
Now $\mtx{T}^{\beta}_{\rm new}$ is a map defined on the
same set of points as all of the neighbors of $\Omega_{\beta}$ and Neumann data can be equated on all sides.

Next, suppose a refined parent does not have a neighbor on one of its sides.  Then on that side we use $\mtx{P}_{\rm down} = \mtx{P}_{\rm up} = \mtx{I}$.  This could happen if the parent is on the boundary of our domain $\Omega$.  For example, suppose $\Omega_{\alpha}$ in Figure \ref{fig:interp_fig} was also refined.  Then $\Omega_{\alpha}$ would not have a neighbor on its Western side.  Additionally, if multiple levels of refinement are done then a refined parent could have no neighbors on one side.  For example, suppose the Northwestern child of $\Omega_{\beta}$ in Figure \ref{fig:interp_fig} was refined.  Then the Northwestern child would not have a neighbor on its Western side since box $\Omega_{\alpha}$ is on a different level of the tree.

Forming the interpolation operators for each side of $\Omega_{\beta}$ before we perform any following merge operations is the easiest approach.  The alternative would be to form an interpolation operator every time two boxes are merged and the nodes do not align.

\subsection{Modifications to solve stage}

On the upwards pass of the solve stage, the fluxes for the particular solution
must be calculated on the same nodes so the particular solution can be
calculated on those nodes.  This is easily achieved by applying the already computed
interpolation operator $\mtx{P}_{\rm down}$ to obtain $\mtx{h}^{\beta}_{\rm new} = \mtx{P}_{\rm down} \mtx{h}^{\beta}_{\rm old}$.


In the downwards pass of the solve stage, the application of the solution operators
results in the approximate solution at the coarse nodes on the Western and Southern sides of $\Omega_{\beta}$.  The solution operator $\mtx{S}^{\beta}_{\rm new}$ now maps the solution on the coarse nodes on the Western and Southern sides of $\Omega_{\beta}$ (and the dense nodes on the Eastern and Northern sides) to the solution on the interior of $\Omega_{\beta}$.  However, we also need the solution on the dense nodes on the Western and Southern sides of $\Omega_{\beta}$.  Let
$\vct{u}_{ge,\rm{new}}$ denote the solution on the boundary of $\Omega_{\beta}$ with
the coarse nodes on the Western and Southern edges.  Then the approximate solution
on the dense nodes is given by $\vct{u}_{ge,\rm{old}} = \mtx{P}_{\rm up}\vct{u}_{ge,\rm{new}}$.

\section{Numerical experiments}
\label{sec:num}

In this section, we present the results of numerical experiments that
illustrate the performance of the scheme proposed.
Section \ref{sec:speed} reports on the computational cost and memory requirements.
Sections \ref{sec:variable}--\ref{sec:tunnel}
report on the accuracy of the proposed solution technique for a
variety of problems where local mesh refinement is required. Finally,
Section \ref{sec:timestep}
illustrates the use of the proposed method in the acceleration of
an implicit time stepping scheme for solving a parabolic partial differential
equation.

For each experiment, the error is calculated by comparing the approximate solution
with a reference solution $\vct{u}_{\rm ref}$ constructed using a highly over resolved grid.
Errors are measured in $\ell^{\infty}$-norm, on all Chebyshev nodes on leaf boundaries.

In all of the experiments, each leaf is discretized using a $p\times p$ tensor product mesh of Chebyshev nodes.
The number of Legendre nodes per leaf edge is set to $q = p-1$.
In all experiments except the one described in Section \ref{sec:tunnel}, the computational domain is
the square $\Omega = [0,1]^{2}$ discretized into $n\times n$ leaf boxes, making the total
number of degrees of freedom roughly $N\approx p^{2}\times n^{2}$
(to be precise, $N = n^{2} \times (p-1)^{2}+2n \times (p-1) + 1$).

The proposed method was implemented in Matlab and all experiments were run on a
laptop computer with a 4 core Intel i7-3632QM CPU running at 2.20 GHz with 12 GB of RAM.

\subsection{Computational speed}
\label{sec:speed}

The experiments in this section illustrate the computational complexity and
memory requirements of the direct solver.  Recall that the asymptotic complexity
of the method scales as nested dissection or multifrontal methods,
with execution times scaling as $O(N^{3/2})$  and $O(N \log N)$ for
the ``build'' and ``solve'' stages, respectively. The asymptotic memory
requirement is $O(N \log N)$.

The computational complexity and memory requirements of the proposed method depend only
on the domain and the computational mesh; the choice of PDE is irrelevant. In the
experiments reported here, we used $\Omega = [0,1]^{2}$ with a uniform mesh.


Figure \ref{Timing_ab} reports the time in seconds for the (a) build and
(b) solve stages of the proposed solution technique when there is a body load (BL),
when the leaf computation is done on the fly as described in Remark \ref{remark:econ} (BL(econ)),
and when there is no body load (NBL). Results for two different orders of discretization ($q = 8$
and $16$) are shown. Notice that as expected the constant scaling factor
for both stages is larger for the higher order discretization.

Figure \ref{Memory} reports on memory requirements. Letting $R$ denote total memory
used, we plot $R/N$ versus the number of discretization points $N$, where $R$
is measured in terms of number of floating point numbers. We see that storing the solution
operators on the leaves is quite costly in terms of memory requirements. The trade-off
to be considered here is whether the main priority is to conserve memory, or to
maximize the speed of the solve stage, cf.~Remark \ref{remark:econ}. As an
illustration, we see that for a problem with $q=8$ and $n=128$, for a total of $10^6$ unknowns,
the solve stage takes 1.7 seconds with the solution operators stored versus 9.8 seconds for
performing the local solves on the fly. In situations where the solution is only desired in
prescribed local regions of the geometry, computing the operators on the fly is ideal.

\begin{remark}
When the underlying BVP that is discretized involves a constant coefficient operator, many
of the leaf solution operators are identical. This observation can be used to greatly
reduce storage requirements while maintaining very high speed in the solve stage.
This potential acceleration was \textit{not} exploited in the numerical experiments
reported.
\end{remark}

\begingroup
\begin{center}
\begin{figure}
\subfigure[]{\includegraphics[height=.45\textwidth]{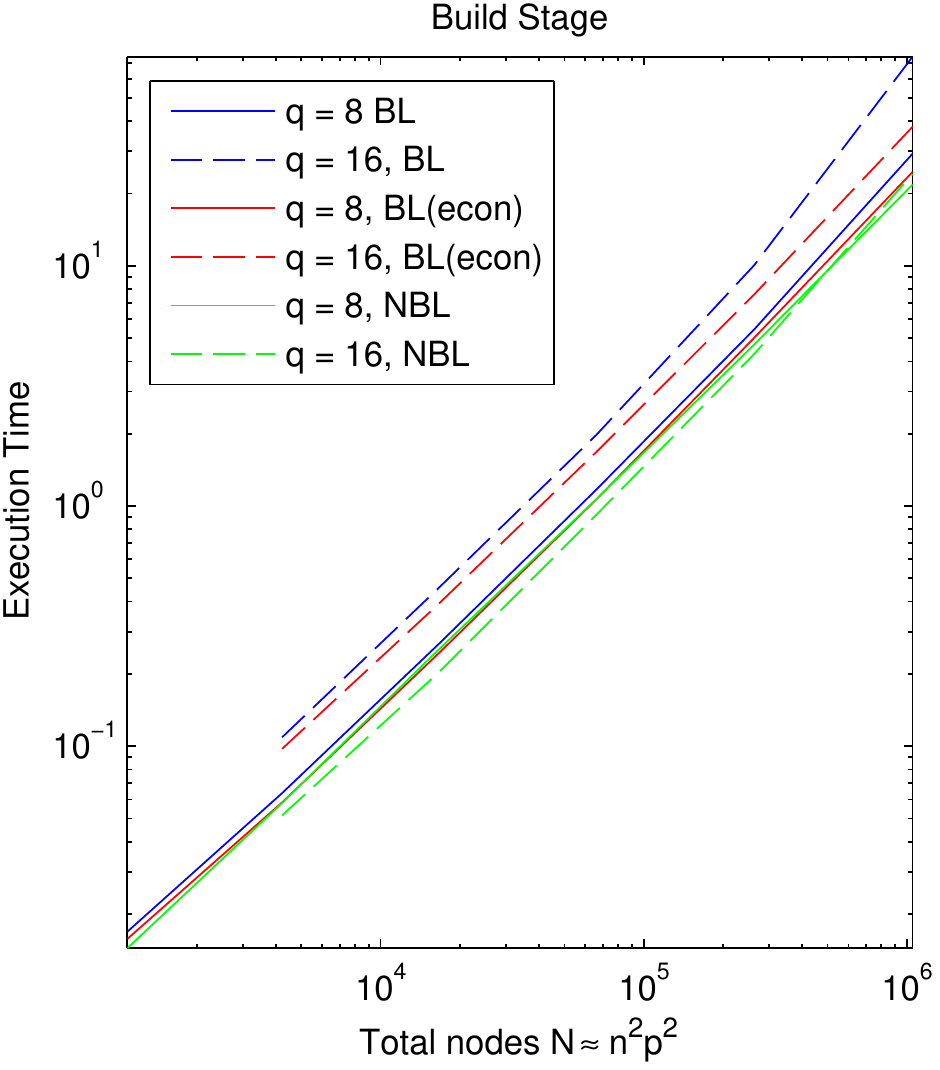}}
\subfigure[] {\includegraphics[height=.45\textwidth]{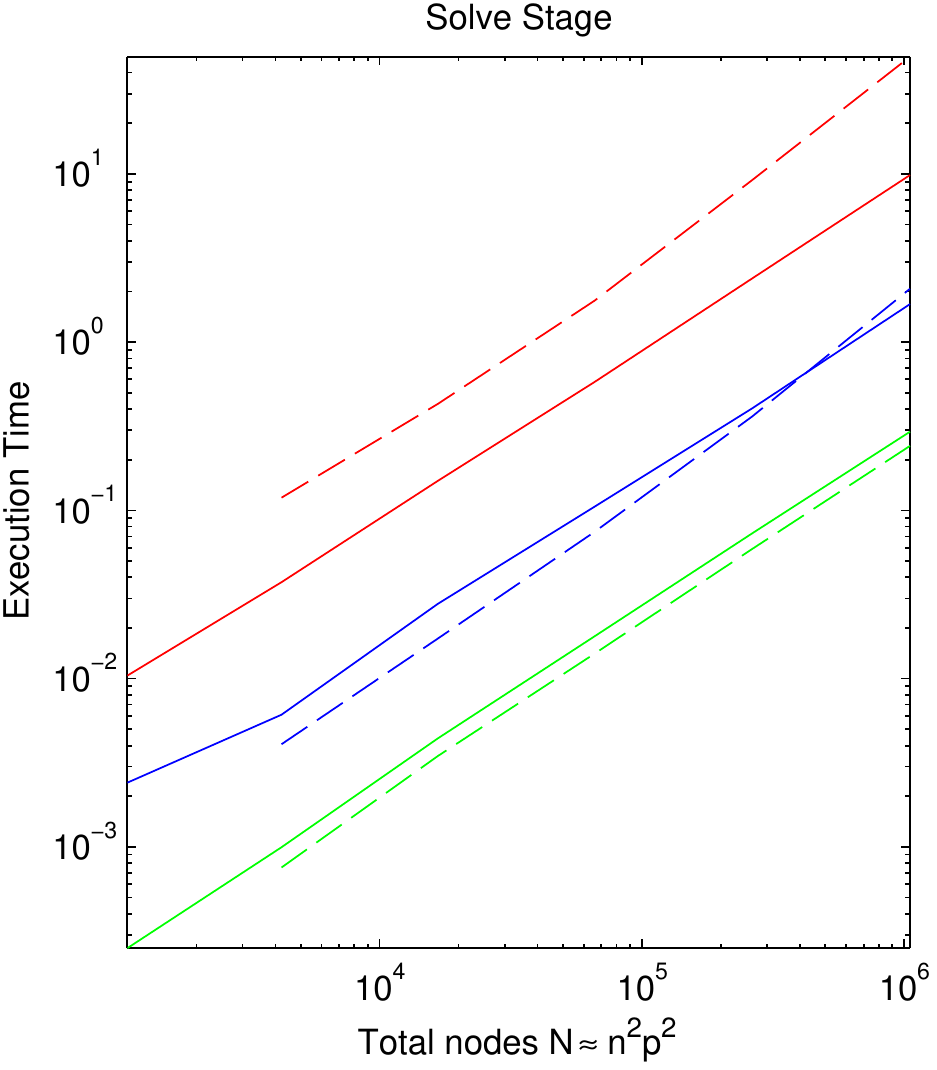}}
\caption{(a) Time to execute build stage for the algorithm with and without a body load.
These algorithms all have complexity $O(N^{1.5})$, and we see that the scaling
factors depend strongly on the order of the method, but only weakly on whether body loads
are included or not. (b) Time to execute the solve stage. Three cases are considered:
NBL is the scheme for problems without a body load.
BL is the scheme for problems with a body load.
BL(econ) is a scheme that allows for body loads, but do not store the
relevant solution operators at the leaves. $p$ denotes the order in
the local Chebyshev grids, and $q=p-1$ is the number of Legendre nodes
on the edge of each leaf.} \label{Timing_ab}
\end{figure}
\end{center}
\endgroup

%


\begin{figure}
\begin{center}
\includegraphics[height=.5\linewidth]{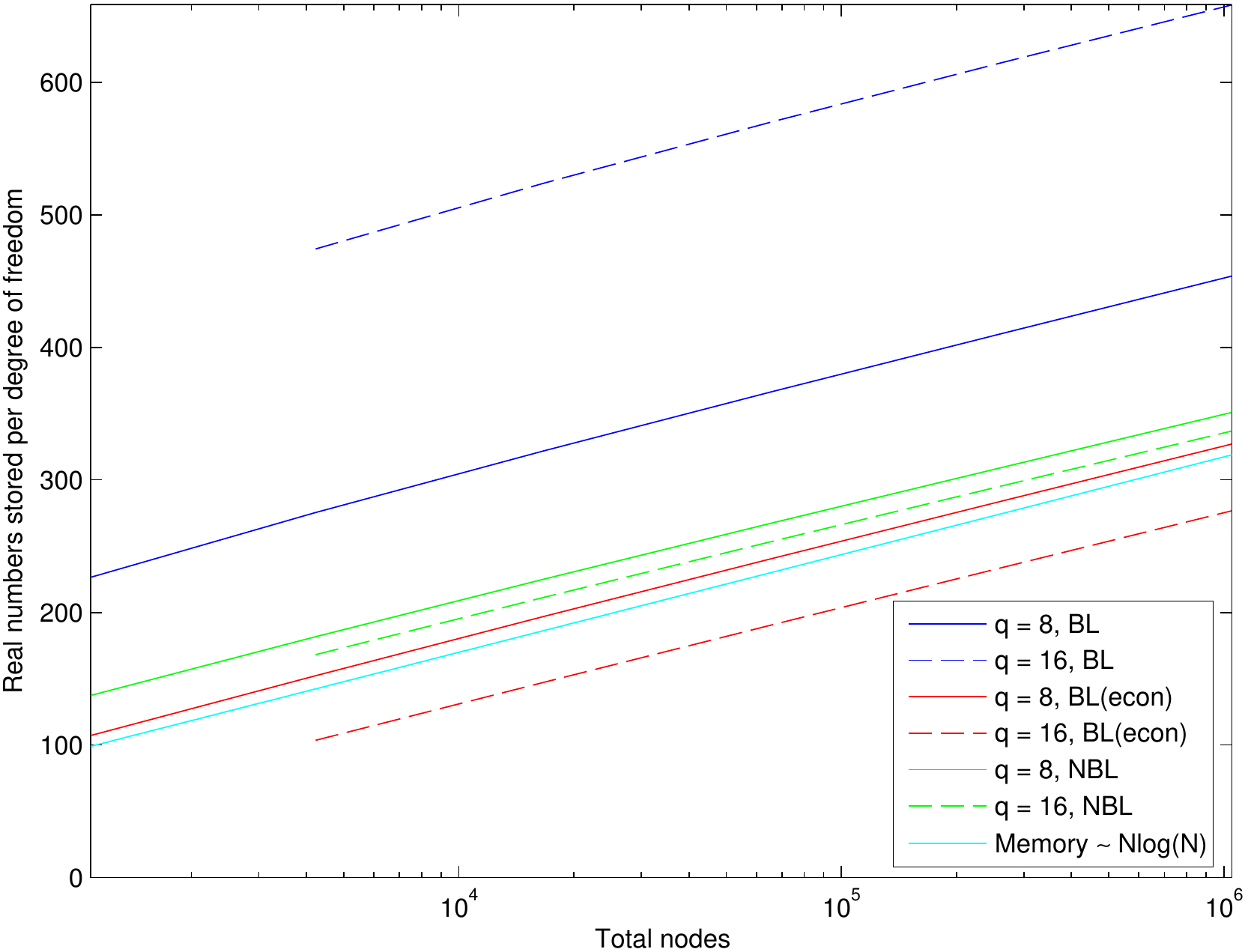}
\caption{Memory requirements. Notation is as in Figure \ref{Timing_ab}.}
\label{Memory}
\end{center}
\end{figure}

\subsection{Variable coefficients}
\label{sec:variable}

In this section, the proposed scheme is applied to the variable coefficient Helmholtz problem
\begin{equation*}
-\Delta u - \kappa^{2}(1-c(\pxx))u = g,\qquad \pxx \in \Omega,
\end{equation*}
where $\Omega = [0,1] \times [0,1]$ and where $c$ is a ``scattering potential.''  The body load is taken to be
a Gaussian given by
$g = \exp(- \alpha |\pxx-\hat {\pxx}|^{2})$ with $\alpha = 300$ and $\hat {\pxx} = [1/4, 3/4]$ while the
variable coefficient is a sum of Gaussians $c(\pxx) = \frac{1}{2}\exp(- \alpha_{2} |\pxx-\hat {\pxx}_{2}|^{2}) +
\frac{1}{2}\exp(- \alpha_{3} |\pxx-\hat {\pxx}_{3}|^{2})$ with $\alpha_{2} = \alpha_{3} =  200$, $\hat {\pxx}_{2} = [7/20, 6/10]$, and $\hat {\pxx}_{3} = [6/10, 9/20]$ for the scattering potential.
We set $\kappa = 40$, making the domain $6.4 \times 6.4$ wavelengths in size.

\begin{figure}
\begin{center}
\includegraphics[height=.5\linewidth]{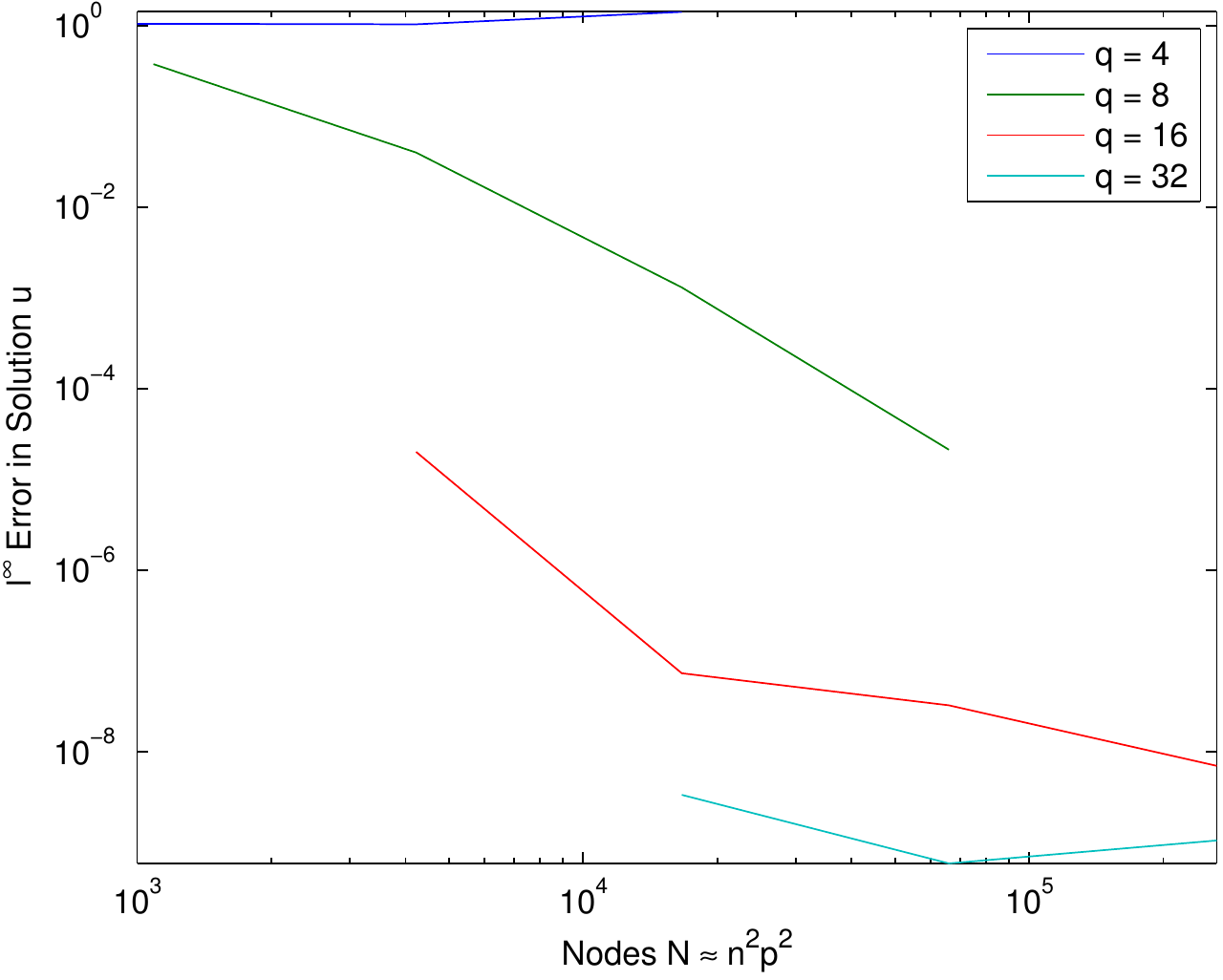}
\caption{The error for the variable coefficient problem described in Section \ref{sec:variable}.
As before, $q$ denotes the number of Legendre nodes along one side of a leaf.}
\label{VarCoef}
\end{center}
\end{figure}

Figure \ref{VarCoef} reports the $l^\infty$ error versus the number of discretization points $N$.
We get no accuracy for $q=4$, but as $q$ is increased, the errors rapidly decrease.


\subsection{Concentrated body load}
In this section, we consider a low frequency ($\kappa = 20$) Helmholtz boundary value problem
\begin{equation*}
-\Delta u - \kappa^{2} u = g,\qquad \pxx \in \Omega,
\end{equation*}
with $\Omega = [0,1] \times [0,1]$ and a very concentrated Gaussian for the body load,
$g = \exp(- \alpha |\pxx-\hat {\pxx}|^{2})$ with $\alpha = 3000$. In this case, we chose
the Dirichlet boundary data to equal the solution to the free space equation $-\Delta u - \kappa^{2} u = g$
with a radiation condition at infinity. In other words, $u$ is the convolution between $g$
and the free space fundamental solution. We computed the boundary data and the reference
solution by numerically evaluating this convolution to very high accuracy.

\begin{figure}
\begin{center}
\includegraphics[height=.5\linewidth]{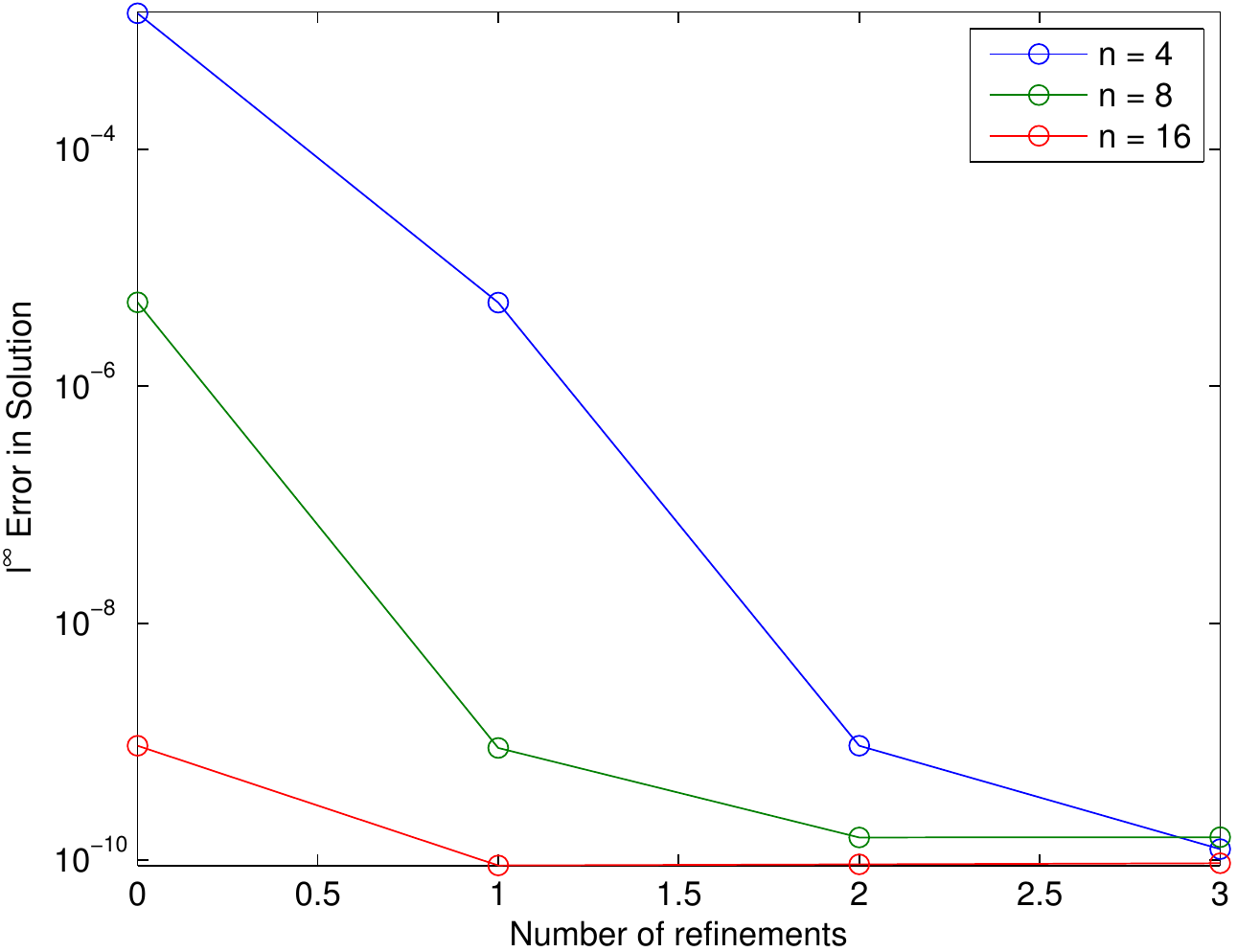}
\caption{Error for Helmholtz equation with $\kappa = 20$ and a very concentrated
body load, demonstrating the ability to improve the solution with refinement.
We use local Chebyshev grids with $17 \times 17$ Chebyshev nodes per leaf,
and $n\times n$ leaves, before refinement. For a problem like this with a concentrated
body load we can improve the error just as much by refining the discretization at the
troublesome location as we can from doubling the number of leaves, which would give the
same grid at the target location.}
\label{Refine}
\end{center}
\end{figure}

To test the refinement strategy, we build a tree first with a uniform grid, i.e. $n\times n$
leaf boxes then add $n_{\rm ref}$ levels of refinement around the point $\hat{\mathbf{x}}$.
Figure \ref{Refine} reports the $l^\infty$ norm of the error versus $n_{\rm ref}$ for
four choices of uniform starting discretization.  When $n = 4$ one level of refinement (i.e. 28 leaf boxes)
results in approximately the same accuracy as when $n = 8$ and no levels of refinement (i.e. 64 leaf
boxes).


\subsection{Discontinuous body load}
\begin{figure}
\begin{center}
\includegraphics[height=.5\linewidth]{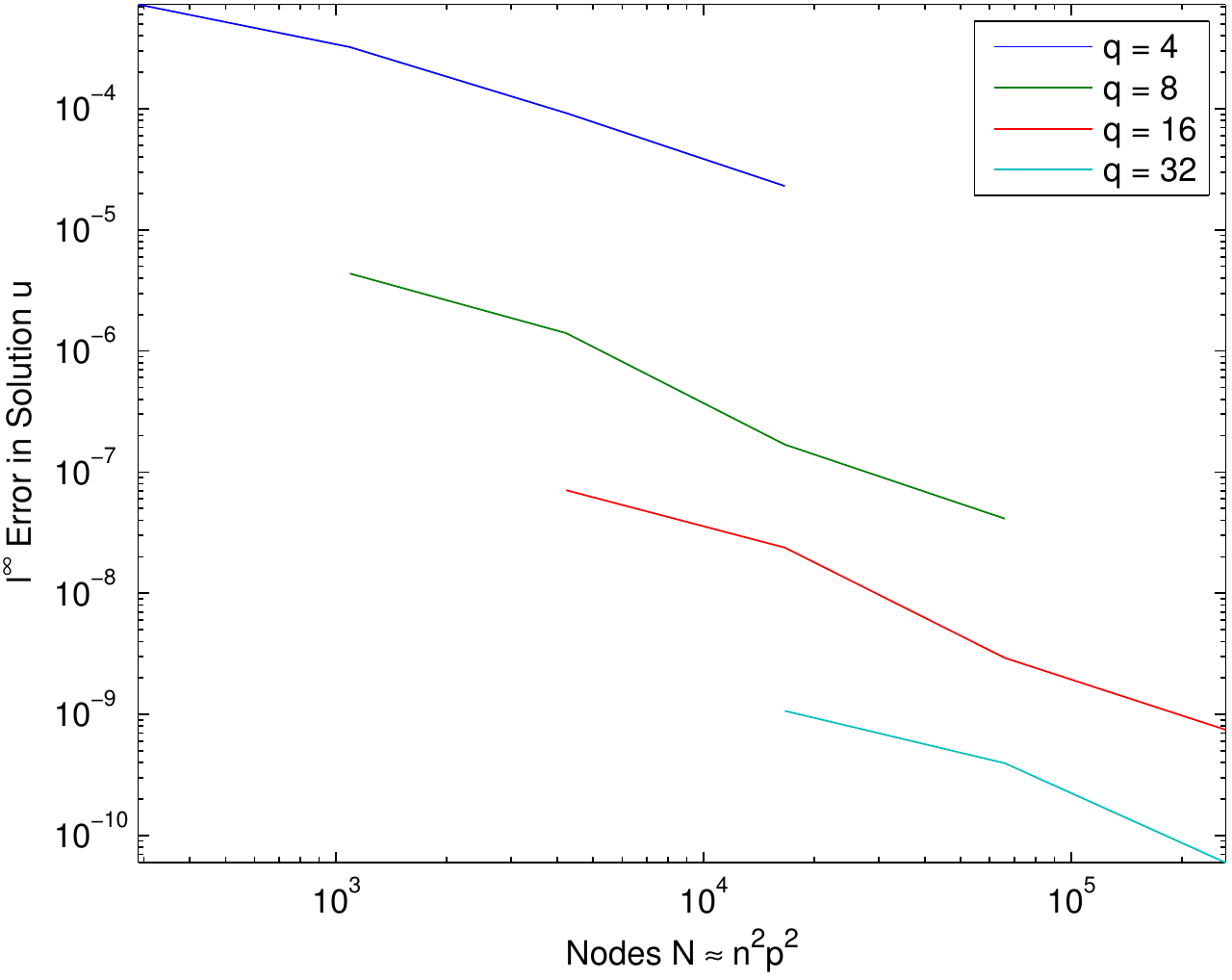}
\caption{The error for a problem with a discontinuous body load.  The discontinuities
align with the edges of the leaves so we still get 10 digits of accuracy.  In the legend,
$q$ denotes the number of Legendre nodes along one side of a leaf.}
\label{Indicator}
\end{center}
\end{figure}

In section, we consider a Poisson boundary value problem on $\Omega = [0,1]^2$
with an indicator function body load $g$ that has support $[1/4, 1/2] \times [1/4, 1/2]$.
Observe that the lines of discontinuity of $g$ coincide with edges of leaves in the
discretization. Figure \ref{Indicator} reports the $l^\infty$ error versus the number
of discretization points $N$ with uniform refinement for four different orders of discretization.
Note that the approximate solution and its first derivative are continuous through the
boundaries of the leaves (even on the boundaries where the jump in the body load occurs)
since the algorithm enforces them by derivation.


\begin{remark}
Applying the scheme to a problem where a discontinuity in the body load does
\textit{not} align with the leaf boundaries results in a very low accuracy approximation
to the solution.  For a problem analogous to the one described in this section, we observed
slow convergence and attained no better than two or three digits of accuracy on the most finely
resolved mesh.
\end{remark}

\subsection{Tunnel}
\label{sec:tunnel}
\begin{figure}
\begin{center}
\includegraphics[height=.5\linewidth]{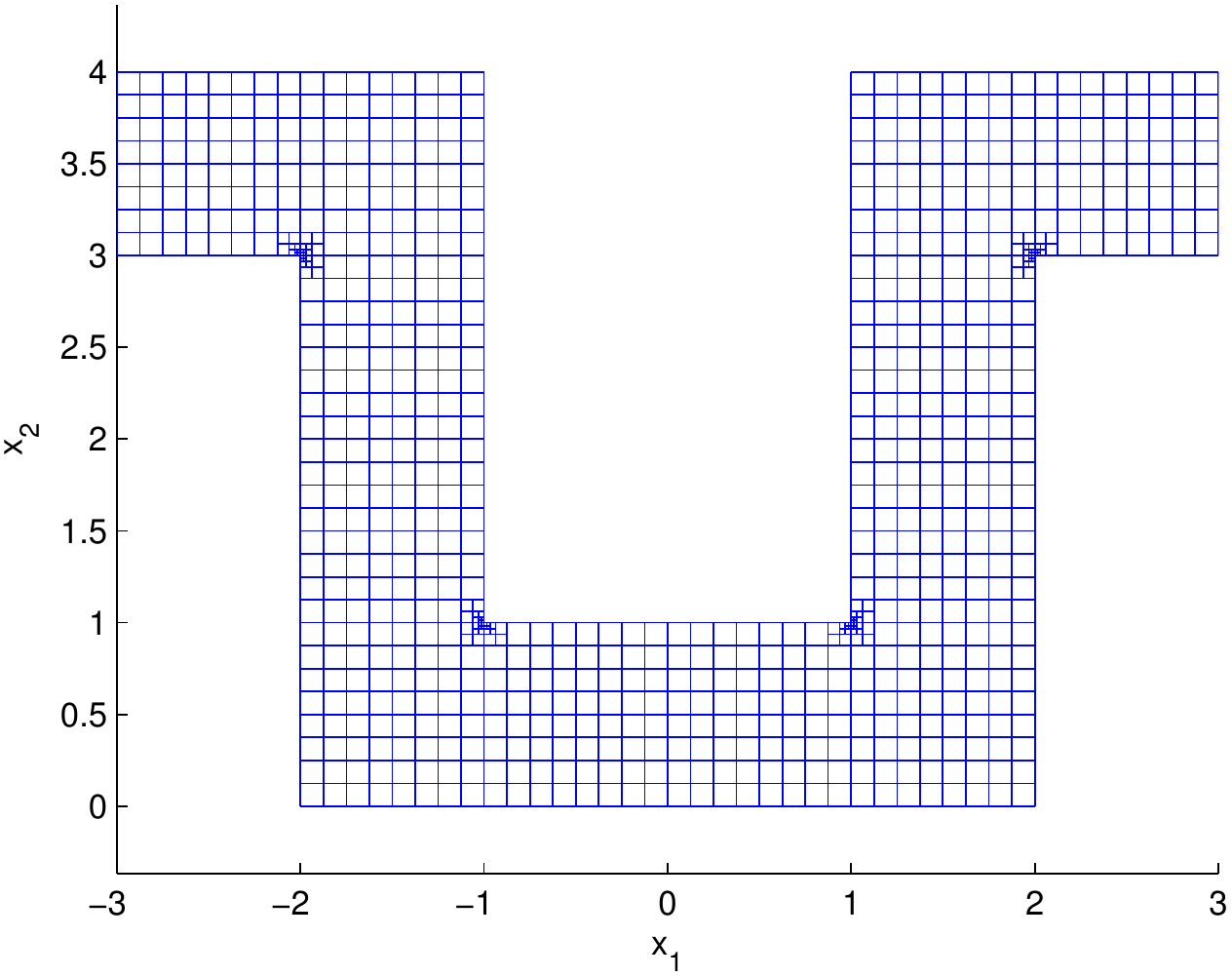}
\caption{The domain used for the tunnel problem.  We solve Helmholtz
equation with $\kappa = 60$, making the tunnel about $10\lambda$ wide
and $115\lambda$ long.  The end caps have fixed Dirichlet data
and the sides of the tunnel have the Dirichlet data set to $u(\pxx) = 0$.}
\label{TGrid}
\end{center}
\end{figure}

This section reports on the performance of the solution technique when applied to
the Helmholtz Dirichlet boundary value problem
\begin{align*}
-\Delta u - \kappa^{2} u &= g,\qquad \pxx \in \Omega,\\
u &= f \qquad \pxx \in \partial \Omega
\end{align*}
with $\kappa = 60$ where the domain $\Omega$ is ``tunnel'' as illustrated in Figure \ref{TGrid}.
The body load is taken to be a Gaussian $g = \exp(- \alpha |\pxx-\hat {\pxx}|^{2})$ with
$\alpha = 300$ and $\hat {\pxx} = [1, 3/4]$. The Dirichlet boundary data is given
by
\begin{equation*}
f(\pxx)= \left\{\begin{array}{cll} 0 & \quad& {\rm for} \ x_1 \neq \pm 3\\
                  \frac{1}{100}\sin(2\pi (x_{2}-3)) & & {\rm for} \ x_1 = 3\\
                  \frac{1}{100}\sin(\pi (x_{2}-3)) & & {\rm for} \ x_1 = -3.
                \end{array}\right.
\end{equation*}
Note that $f(\pxx)$ is continuous on $\partial \Omega$ and with this choice of wave number $\kappa$
the domain $\Omega$ is about 10 wavelengths wide and 115 wavelengths long.  The presence of the re-entrant corners
results in a solution that has strong singularities which require local refinement in order
for the method to achieve high accuracy.

\begin{figure}
\begin{center}
\includegraphics[height=.5\linewidth]{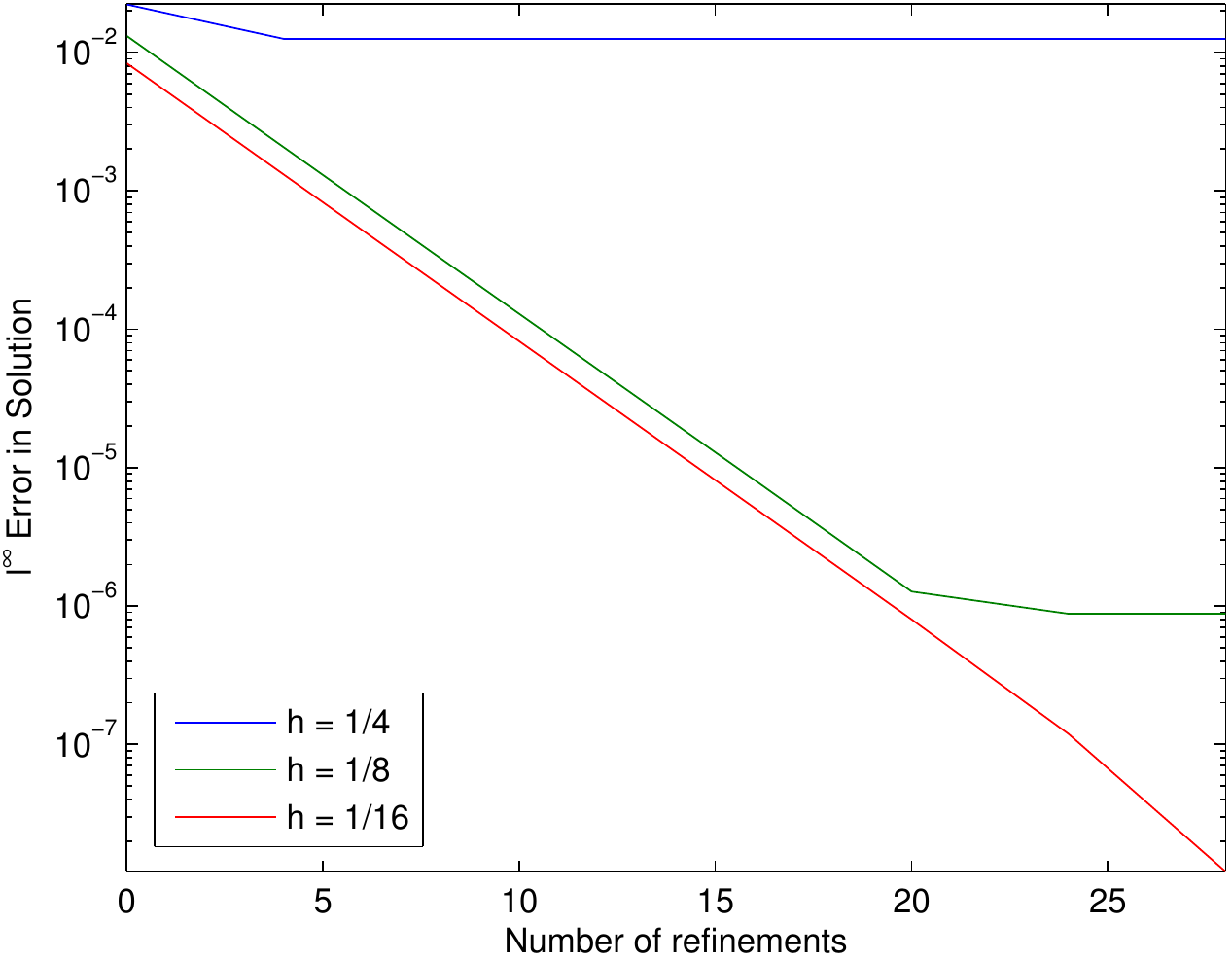}
\caption{Error for Helmholtz equation with $\kappa = 60$ on the tunnel.  The end caps have fixed Dirichlet data and the sides of the tunnel have the Dirichlet data set to $f(\pxx) = 0$.  A Gaussian $g = \exp(- \alpha |\pxx-\hat {\pxx}|^{2})$ with $\alpha = 300$ located at $\hat {\pxx} = [1, 3/4]$ is used for the body load.}
\label{TError}
\end{center}
\end{figure}

Figure \ref{TError} reports the $l^\infty$ error versus the number of refinements into the
corners with three choices of coarse grid.  We use $q = 16$ for all examples and $h$ gives the width and height of each leaf box.  When $h= 1/4$, the discretization is only sufficient to resolve the Helmholtz equation with $\kappa = 60$ within 1\% of the exact solution. 
When $h = 1/8$, the solution
technique stalls at 5 digits of accuracy independent of the number of refinement levels.

\begin{remark}[Symmetries]
This problem is rich in symmetries that can be used to accelerate the build stage.
In our implementation, we chose to exploit the fact that the tunnel is made
up of four L-shaped pieces glued together.   The DtN operator and corresponding solution
operators were constructed for one L-shape.  Then creating the solver for the entire geometry
involved simply gluing the 4 L-shaped geometries together via \textit{three} merge operations.
\end{remark}

\subsection{A parabolic problem}
\label{sec:timestep}
Our final numerical example involves a convection-diffusion initial value problem on $\Omega = [0,1]^2$ given by
\begin{align*}
\Big( \epsilon \Delta - \frac{\partial}{\partial x_{1}} \Big) u(\pxx,t) &= \frac{\partial u}{\partial t}, &\qquad\pxx\in \Omega, \ t>0\\
u(\pxx,0) &= \exp(- \alpha |\pxx-\hat {\pxx}|^{2}),&\qquad\pxx\in \Omega.
\end{align*}
We imposed zero Neumann boundary conditions on the south and north boundaries ($x_{2} = 0,1$) and
periodic boundary conditions on the west and east boundaries ($x_{1} = 0,1$).
These boundary conditions correspond to fluid flowing through a periodic channel where no fluid can
exit the top or bottom of the channel.
To have a convection dominated problem, we chose $\epsilon = 1/200$.
Finally, the parameters in the body load were chosen to be $\alpha = 50$ and $\hat {\pxx} = [1/4, 1/4]$.


Applying the Crank-Nicolson time stepping scheme with a time step size $k$ results in having to solve the
following elliptic problem at each time step:
\begin{equation}
\label{eq:CN}
\Big ( \frac{1}{k} I - \frac{1}{2} A \Big) u_{n+1} = \Big ( \frac{1}{k} I + \frac{1}{2} A \Big) u_{n},
\end{equation}
where $A = \epsilon \Delta - \partial/\partial x_{1}$ is our partial differential operator.

Observe that the algorithm does not change for this problem.  The build stage execution time and memory requirement are identical to those seen in Section \ref{sec:speed}.  The execution time for the solve stage for each individual time step is nearly identical to the solve stage execution time shown in Section \ref{sec:speed}.  The only new step in the solve stage is the need to evaluate $(I/k + A/2)u_{n}$ at each time step.

\begin{figure}
\begin{center}
\includegraphics[height=.5\linewidth]{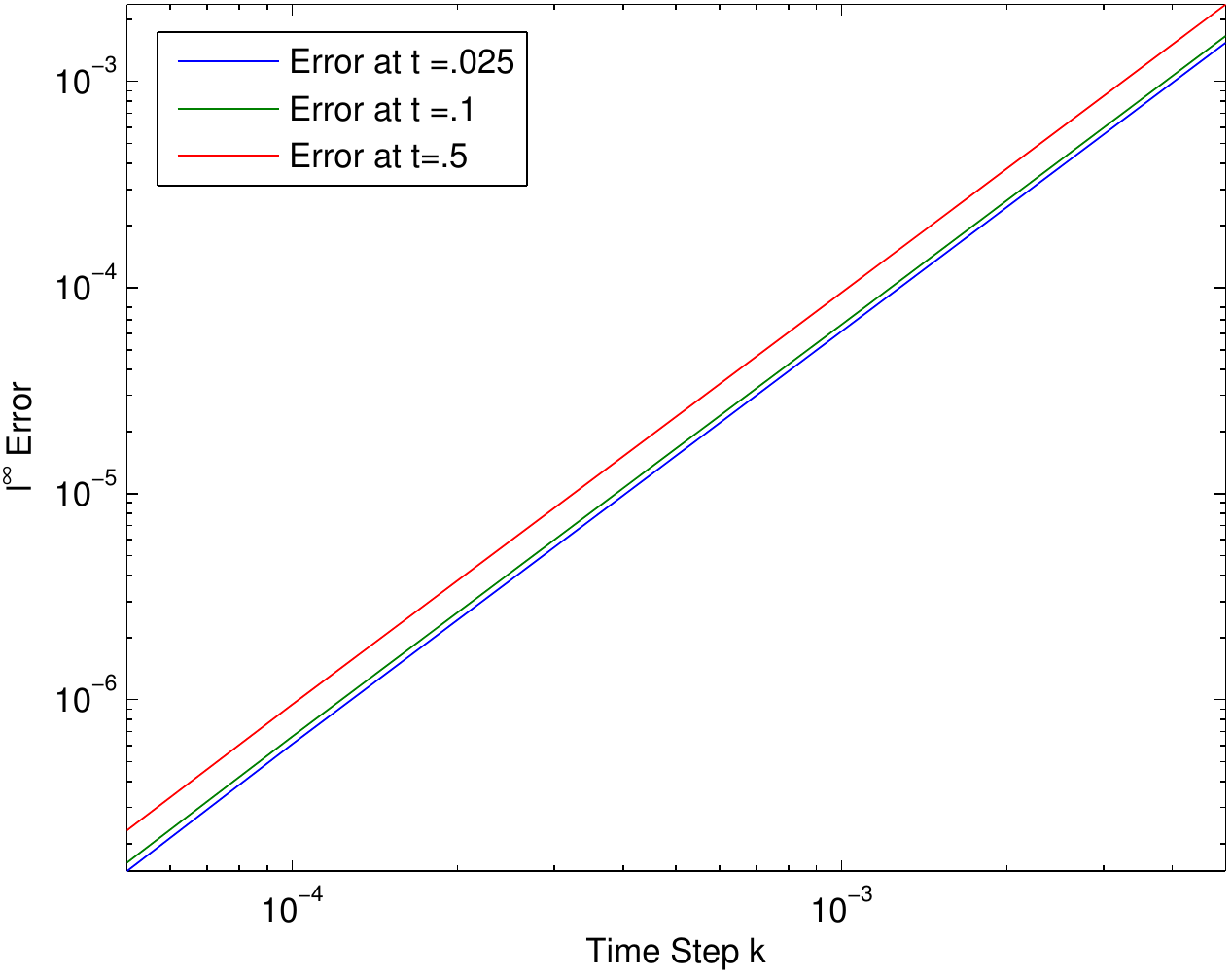}
\caption{Error for the convection-diffusion equation described in Section \ref{sec:timestep}.
The error is estimated by comparing against a highly over-resolved solution.}
\label{ConDiff}
\end{center}
\end{figure}

Figure \ref{ConDiff} reports the $l^\infty$ error vs.~the time step size $k$ at three different
times $t = 0.025,\ 0.1,$ and $0.5$.   Note that even with a low order time stepping
scheme, it is still feasible to high accuracy (i.e. use small time steps)
since processing each time step is very inexpensive.

We use $16$ leaf boxes per side with $q=16$.  This gives more than enough nodes to obtain the accuracy shown in Figure \ref{ConDiff} and the error shown is limited by the accuracy of Crank-Nicolson and not by the discretization in space.


\section{Concluding remarks}
We have described an algorithm for solving non-homogenous linear elliptic PDEs in two
dimensions based on a multidomain spectral collocation discretizations. The solver is
designed explicitly for being combined with a nested dissection type direct solver. Its
primary advantage over existing methods is that it enables the use of very high order
local discretization without jeopardizing computational efficiency in the direct solver.
The scheme is an evolution on previously published methods
\cite{2012_spectralcomposite,2013_martinsson_ItI,2013_martinsson_DtN_linearcomplexity}.
The novelty in this work is that the scheme has been extended to allow for problems
involving body loads, and for local refinement.

The scheme is particularly well suited for executing the elliptic solve required solving
parabolic problems using implicit time-stepping techniques in situations where the domain
is fixed, so that the elliptic solve is the same in every time step. In this environment,
the cost of computing an explicit solution operator is amortized over many time-steps, and
can also be recycled when the same equation is solved for different initial conditions.

The fact that the method can with ease incorporate high order local discretizations,
and allows for very efficient implicit time-stepping appears to make it particularly
well suited for solving the Navier-Stokes equations at low Reynolds numbers. Such a
solver is currently under development and will be reported in future publications.
Other extensions currently under way includes the development of adaptive refinement
criteria (as opposed to the supervised adaptivity used in this work), and the extension
to problems in three dimensions, analogous to the work in \cite{2016_martinsson_HPS_3D}
for homogeneous equations.

\vspace{7mm}

\noindent
\textit{\textbf{Acknowledgements:}} The research reported was supported by DARPA,
under the contract N66001-13-1-4050, and by the NSF, under the contracts DMS-1407340
and DMS-1620472.

\bibliography{main_bib}
\bibliographystyle{amsplain}

\end{document}